\documentclass[12pt]{amsart} 
\usepackage{amsmath,amstext,amsfonts,amssymb} 

\newtheorem{theorem}{\bf Theorem}[section] 
 
\newtheorem{corol}[theorem]{\bf Corollary} 
\newtheorem{defi}[theorem]{\bf Definition} 
\newtheorem{example}[theorem]{\bf Example} 
\newtheorem{lemma}[theorem]{\bf Lemma} 
\newtheorem{notation}[theorem]{\bf Notation} 
\newtheorem{prop}[theorem]{\bf Proposition} 
\newtheorem{rem}[theorem]{\bf Remark} 
\large 

\newcommand{\abs}[1]{\lvert#1\rvert} 
\newcommand{\ctl}{\centerline} 
\newcommand{\ifff}{{if and only if }} 
\newcommand{\nd}{{\text{ and }}} 
\newcommand{\noi}{{\noindent}} 
\newcommand{\bgop}{{\bigoplus}} 
\newcommand{\Equi}{\Longleftrightarrow} 
\newcommand{\la}{{\langle}} \newcommand{\ra}{{\rangle}} 
\newcommand{\lf}{{\lfloor}} \newcommand{\rf}{{\rfloor}} 
\newcommand{\pfi}{{perf-irre\-du\-cible}}

\DeclareMathOperator{\ord}{ord} 
\DeclareMathOperator{\Oo}{O} 
\DeclareMathOperator{\perf}{perf\,rk} 
\DeclareMathOperator{\rk}{rk} 
\DeclareMathOperator{\sgn}{sgn} 
\DeclareMathOperator{\supp}{supp} 
\newcommand{\ta}{{{}^t\!}} 

\renewcommand{\a}{{\alpha}} 
\newcommand{\Lb}{{\Lambda}} 
\newcommand{\lb}{{\lambda}} 
\newcommand{\s}{{\sigma}} 
\newcommand{\va}{{\varepsilon}} 

\newcommand{\D}{{\mathbb D}} 
\newcommand{\E}{{\mathbb E}} 
\newcommand{\R}{{\mathbb R}} 
\newcommand{\Z}{{\mathbb Z}} 

\newcommand{\cB}{{\mathcal B}} 
\newcommand{\cC}{{\mathcal C}} 
\newcommand{\cE}{{\mathcal E}} 
\newcommand{\cL}{{\mathcal L}} 
\newcommand{\cP}{{\mathcal P}} 
\newcommand{\cRR}{{\mathcal R}} 

\DeclareMathOperator{\Aut}{Aut} 
\DeclareMathOperator{\End}{End} 
\DeclareMathOperator{\Id}{Id} 
\DeclareMathOperator{\Mat}{Mat} 
\DeclareMathOperator{\Sym}{Sym}

\begin{document}

\title{On Perfection Relations in Lattices} 
\author[A.-M. Berg\'e, J. Martinet] 
{Anne-Marie Berg\'e and Jacques Martinet ($^*$)} 
\keywords{Euclidean lattices, Quadratic forms, kissing number, 
\newline 
(*) Laboratoire A2X (UMR\,5465, FR\,2254), 
C.N.R.S et Universit\'e Bordeaux~1} 
\begin{abstract} 
Let $\Lb$ be a lattice in a Euclidean space $E$, 
with kissing number~$s$ and perfection rank~$r$, 
that is, the rank in $\End^{\text{sym}}(E)$ 
of the set of orthogonal projections to minimal vectors of~$\Lb$. 
This defines a space of \emph{perfection relations}, 
of dimension $s-r$. We focus on ``short relations'',  
in connection with the index theory, previously developed 
by Watson, Ry\v skov, Zahareva and the second author 
in \cite{W}, \cite{R}, \cite{Z} and \cite{M1}. 
\end{abstract} 
\maketitle

\section{Introduction}\label{secintro} 

Let $(E,x\cdot y)$ be a Euclidean space, of dimension~$n$. 
For every subspace 
$F$ of $E$, denote by $p_F$ the orthogonal projection 
to~$F$. Given a set $\cL$ of $s$~lines in~$E$, 
a \emph{perfection relation on $\cL$} is a relation 
(in the set $\End^s(E)$ of symmetric endomorphisms of $E$) 
$\sum_{L\in \cL}\,\a_L\,p_L=0$ with real coefficients~$\a_L$. 
In practice, we consider the set $S$ of norm~$1$ vectors $\pm x$ 
which belong to the lines of $\cL$, 
and set $N(x)=x\cdot x$ and $p_x=p_L$. 
Since $N(x)=1$, we then have $p_x(y)=(x\cdot y)\,x$ for every $y\in E$. 
The \emph{perfection rank of $\cL$} is the rank $r=\perf\cL$ 
in $\End^s(E)$ of the set $p_L$, $L\in\cL$. 
We say that this family is \emph{perfect} if $r=\frac{n(n+1)}2$. 

\smallskip 

In the forthcoming sections, we shall apply the definitions above 
to the set $S=S(\Lb)$ of minimal vectors of a lattice $\Lb$ in~$E$. 
In this case, we recover the notion of a \emph{perfect lattice}. 
Without loss of generality, we may assume that $\Lb$ is generated 
by those of its minimal vectors which are involved 
in the perfection relation. In particular, $\Lb$ is 
\emph{well   rounded},  i.e., we have $\rk S=\dim\Lb$. 

The set of possible structures for $\Lb/\Lb'$ where $\Lb'$ is 
generated by $n$~independent minimal vectors of $\Lb$ will 
play a major r\^ole in this paper, and in particular, 
the maximal value $\imath$ of the index $[\Lb:\Lb']$. 
It turns out that any perfection relation may be written 
in the form 
$$\sum_{i=1}^m\lb_i\,p_{e_i}=\sum_{j=1}^{m'}\lb'_j\,p_{e'_j}$$ 
where both the  systems $\{e_i\},\{e'_j\}$ are of rank~$n$
(we then denote by $\Lb_0,\Lb'_0$ the lattices they generate) 
and the coefficients $\lb_i,\lb'_i$ are strictly positive. 
We shall focus on the simplest case when $m=m'=n$, 
but even in this simple case, we shall obtain complete 
classification results only under one of the assumptions 
``$[\Lb:\Lb_0]\le 4$'' or ``$\Lb/\Lb_0$ is $2$-elementary'', 
which however covers all dimensions $n\le 7$. 

\smallskip 

Section~\ref{seceuclid} is devoted to the proof of technical results 
on perfection relations in Euclidean spaces 
and Section~\ref{seclatperfrel} to the particular case of lattices. 
In Sections~\ref{seclat2elem}, \ref{secindex3} and~\ref{secindex4}, 
we classify lattices for which $\Lb/\Lb_0$ is $2$-elementary 
or cyclic of order $3$ or~$4$. We discuss various complements 
(action of groups, dimension~$8$,\,...) in Section~\ref{seccompl}. 

\smallskip 

\noi{\sc Acknowledgements}. We would like to thank the authors 
of the {\em PARI} system, and more specially Christian Batut 
and Karim Belabas for their help in applying {\em PARI} 
to lattices.

\section{Perfection Relations in Euclidean Spaces.}\label{seceuclid} 

In this section, we consider perfection relations on a set $\cL$ of lines 
(or on a symmetric set of vectors of norm~$1$). 
Except in the last assertion of Lemma~\ref{lemAMB2}, 
we do not make use of lattices. 

A perfection relation $\sum_{x\in S/\pm}\,\lb_x\,p_x=0$ may be written 
$$\sum_{x\in T/\pm}\lb_x\,p_x=\sum_{x\in T'/\pm}\lb'_x\,p_x$$ 
with \emph{strictly positive} coefficients $\lb_x,\,\lb'_x$. 

\begin{lemma}\label{lemrank} 
With the notation above, $T$ and $T'$ span the same subspace of~$E$. 
In particular, they have the same rank. 
\end{lemma} 

\begin{proof} 
Let $F$ be the span of $T$ and $F'$ that of $T'$. 
For every $y\in {F'}^\perp$, we have 
$\big(\sum_{x\in T/\pm}\lb_x\,p_x(y)\big)\cdot y=0$, i.e. 
$\sum_{x\in T/\pm}\lb_x\,(x\cdot y)^2=0$, 
which implies $x\cdot y=0$ for all $x\in T$, hence $x\in F^\perp$. 
We thus have ${F'}^\perp\subset F^\perp$, and similarly 
$F^\perp\subset{F'}^\perp$, hence $F^\perp={F'}^\perp$, 
i.e. $F=F'$. 
\end{proof} 

\begin{rem}\label{remsigns} 
{\small\rm 
Let $(e_1,\dots,e_n)$ be a (unitary) basis for $E$. 
Set $u=\sum_i\,\lb_i\,p_{e_i}$, with $\lb_1,\dots,\lb_n\in\R$. 
Then Sylvester's law of inertia applied to the quadratic form 
$u(x)\cdot x=\sum\,\lb_i\,(e_i\cdot x)^2$ 
shows that the numbers of $\lb_i$ which are $>0$, $<0$ or zero 
depend only on~$u$.} 
\end{rem} 

It results from Lemma~\ref{lemrank} that we may restrict ourselves 
to perfection relations in which both $T$ and $T'$ span~$E$. 
Then such a relation involves at least $2n$ lines, 
and it is easy to check that when exactly $2n$~lines are involved, 
this is then unique up to proportionality except if it comes 
from two relations in two strict subspaces of~$E$. 
Most of the time, we shall assume that no such subspaces exist. 
Then $\perf\,(T\cup T')=2n-1$. 
Proposition~\ref{propreduc} below describes a situation 
in which perfection relations on two complementary spaces occur. 

\begin{prop}\label{propreduc} 
Let $\cB=(e_1,\dots,e_n)$ and $\cB'=(e'_1,\dots,e'_n)$ 
be two bases for $E$ and let $\lb_1,\dots,\lb_n,\lb'_1,\dots,\lb'_n$ 
be strictly positive real numbers 
such that $\sum_{i=1}^n\,\lb_i\,p_{e_i}=\sum_{i=1}^n\,\lb'_i\,p_{e'_i}$. 
Assume that there exists two partitions 
$\{1,\dots,n\}=I_1\cup I_2=J_1\cup J_2$ 
such that, for $k=1$ or $2$, each $e'_i,\,i\in J_k$ belongs 
to the span $E_k$ of $\{e_i,\,i\in I_k$\}. 
Then $J_k$ and $I_k$ have the same cardinality, 
and we have the two perfection relations 
$$\sum_{i\in I_1}\lb_i\,p_{e_i}=\sum_{i\in J_1}\lb'_i\,p_{e'_i}\ \nd\ 
\sum_{i\in I_2}\lb_i\,p_{e_i}=\sum_{i\in J_2}\lb'_i\,p_{e'_i}\,.$$ 
\end{prop} 

\begin{proof} 
Let $u$ be the symmetric endomorphism defined by either side 
of the equality 
$$\sum_{i\in I_1}\lb_i\,p_{e_i}-\sum_{i\in J_1}\lb'_i\,p_{e'_i}= 
\sum_{i\in J_2}\lb'_i\,p_{e'_i}-\sum_{i\in I_2}\lb_i\,p_{e_i}\,.$$ 
For all $x\in E$, we have $u(x)\in E_1\cap E_2=\{0\}$, 
i.e. $u$ is zero. 
\linebreak 
By Lemma~\ref{lemrank}, 
$\rk\{e_i,\,i\in I_k\}=\rk\{e'_i,\,i\in J_k\}$, 
hence $\abs{I_k}=\abs{J_k}$ for $k=1,\,2$. 
\end{proof} 

\begin{defi}\label{defirred} 
We say that the set $\{e_i,e'_j\}$ is \emph{\pfi} 
if no such system of partitions exists. 
\end{defi} 

Returning to the previous notation , we now prove a characterization 
of perfection relations involving two bases for~$E$. 
Recall that given a basis $\cB=(e_1,\dots,e_n)$ 
for $E$ with dual basis $\cB^*=(e_1^*,\dots,e_n^*)$ 
(i.e., $e_i\cdot e_j^*=\delta_{i,j}$), for every $x\in E$, 
the scalar products $x\cdot e_i^*$ are the components of $x$ on the $e_i$. 

\begin{lemma}\label{lembases} 
Let $\cB=(e_1,\dots,e_n)$ and $\cB'=(e'_1,\dots,e'_n)$ 
be two bases for $E$ and let $\lb_1,\dots,\lb_n,\lb'_1,\dots,\lb'_n$ 
be real numbers. Then the following conditions are equivalent: 
\begin{enumerate} 
\item 
$\sum_{i=1}^n\,\lb_i\,p_{e_i}=\sum_{i=1}^n\,\lb'_i\,p_{e'_i}$. 
\item 
$\forall\,j$, $\lb_j\,e_j=\sum_{i=1}^n\,\lb'_i(e'_i\cdot e_j^*)e'_i$. 
\item 
$\forall\,k$, $\lb'_k\,e'_k=\sum_{i=1}^n\,\lb_i(e_i\cdot {e'}_k^*)e_i$. 
\item 
$\forall\,j,\,\forall\,k$, 
$\lb_j(e_j\cdot {e'}_k^*)=\lb'_k(e'_k\cdot e_j^*)$. 
\end{enumerate}\end{lemma} 

\begin{proof} 
Both sides of (1) are endomorphisms of $E$, which are equal \ifff 
they coincide on some basis. Taking the values of both sides on $\cB^*$ 
(resp. ${\cB'}^*$) gives $(1)\Leftrightarrow(2)$ 
(resp. $(1)\Leftrightarrow(3)$). 
Then we observe that (2) is a collection of equalities between 
$n$~pairs of vectors of $E$, and the vectors of both sides are equal 
\ifff they have the same scalar products with the $n$~vectors of some 
basis. Using the basis ${\cB'}^*$, we obtain the equivalence 
of the $n$~equalities in~(2) and the $n^2$~equalities in~(4). 
\end{proof} 

\begin{lemma}\label{lemAMB1} 
Suppose that the $\lb_i$ and $\lb'_i$ satisfy the equivalent conditions 
of Lemma~\ref{lembases}. Let 
$$A_i=1-\sum_{k=1}^n\,(e_i\cdot {e'}_k^*)^2\ \nd\ 
A'_i=1-\sum_{j=1}^n\,(e'_i\cdot e_j^*)^2\,.$$ 
Then 

\ctl{$\sum_{i=1}^n\,\lb_i=\sum_{i=1}^n\,\lb'_i$} 

\noi and 

\ctl{$\sum_{i=1}^n\,\lb_i\,A_i=\sum_{i=1}^n\,\lb'_i\,A'_i=0\,.$} 
\end{lemma} 

\begin{proof} 
The first assertion results from the fact 
that projections to a line have trace~$1$. 

Taking the scalar product with $e_j^*$ 
of the two sides of formula (2) in Lemma~\ref{lembases} 
yields the formula 
$$\lb_j=\sum_{i=1}^n\,\lb'_i\,(e'_i\cdot e_j^*)^2\,.$$ 
By summation on $j$, the left hand side becomes $\sum_j\,\lb_j$, 
equal to $\sum_i\,\lb'_i$ by the assertion above, 
which proves that $\sum_{i=1}^n\,\lb'_i\,A'_i=0$. 
Exchanging the systems $(\lb_i,e_i)$ and $(\lb'_i,e'_i)$ 
completes the proof. 
\end{proof} 

\begin{lemma}\label{lemAMB2} 
Suppose moreover that the $\lb_i$ and $\lb'_i$ are strictly positive. 
Then: 
\begin{enumerate} 
\item 
$\forall\,j,\,\forall\,k$, $e_j\cdot {e'}_k^*$ and $e'_k\cdot e_j^*$ 
have the same sign or are both zero. 
\item 
$\forall j$, 
$\sum_{i=1}^n\,\abs{e_j\cdot {e'}_i^*}\,\abs{e'_i\cdot e_j^*}=1$. 
\item 
$\forall i$, 
$\sum_{j=1}^n\,\abs{e_j\cdot {e'}_i^*}\,\abs{e'_i\cdot e_j^*}=1$. 
\item Assume that the $(e_i)$ and the $(e'_j)$ 
are minimal vectors in some lattice. 
Then for all $j$ and all $k,\,h$ such that 
$e_j\cdot{e'_k}^*\ne 0$ and $e_j\cdot{e'_h}^*\ne 0$, we have 
$\frac{\lb'_k}{\lb'_h}\,\frac{\abs{e'_k\cdot e_j^*}} 
{\abs{e'_h\cdot e_j^*}}= 
\frac{\abs{e_j\cdot{e'_k}^*}}{\abs{e_j\cdot{e'_h}^*}} 
\ge \frac 1{\a_n}$, 
where $\a_n$ is the maximal value for the annihilator 
of $\Lb/\Lb'$ for well-rounded 
\linebreak 
$n$-dimensional lattices 
having the same minimum. 
\newline 
{\small\rm 
[We have $\a_n=1,2,3,4,6$ for $n\le 3$, $n=4\text{\,or\,}5$,
$n=6$, $n=7$ $,n=8$ respectively; see \cite{M1}, Table~11.1.]} 
\end{enumerate}\end{lemma} 

\begin{proof} 
The first assertion results from the last assertion 
of Lemma~\ref{lembases}. 
Replacing $\lb'_i\,(e'_i\cdot e_j^*)$ by $\lb_j\,(e_j\cdot {e'}_i^*)$ 
in the displayed formula which occurs in the proof of Lemma~\ref{lemAMB1}, 
dividing both sides by $\lb_j$ and using (1), we obtain (2). 
Exchanging the systems $(\lb_i,e_i)$ and $(\lb'_i,e'_i)$ 
proves~(3). By \ref{lembases},\,(4), 
$\frac{\lb'_k}{\lb'_h}\,\frac{\abs{e'_k\cdot e_j^*}} 
{\abs{e'_h\cdot e_j^*}}= 
\frac{\abs{e_j\cdot{e'_k}^*}}{\abs{e_j\cdot{e'_h}^*}}$. 
Write $e_j=\frac 1d\,\sum_{i=1}^n\,a_i e'_i$ with coprime integers 
$a_i,d$. Each $e'_i$ with $a_i\ne 0$ may be written as a combination 
of $e_j$ and the $e'_\ell,\ell\ne i$ with denominator $a_i$, 
so that $\abs{a_i}\le\a_n$, which shows that the ratios 
$\frac{\abs{a_k}}{\abs{a_h}}$ ($a_k,a_h\ne 0$) are bounded from below 
by $\frac 1{\a_n}$. 
\end{proof} 

\begin{lemma}\label{lemAMB3} 
With the notation and hypotheses of Lemma~\ref{lemAMB2}, 
let 
\linebreak 
$j$ and $k$ be two indices such that, for all~$i$, 
$e'_k\cdot e_i^*$ is non-zero. Then 
$$\sum_i\,\abs{e'_k\cdot e_i^*}\,\abs{e_i\cdot{e'_k}^*}\, 
\Big(\frac{\lb'_j}{\lb'_k}\,
\frac{(e'_j\cdot e_i^*)^2}{(e'_k\cdot e_i^*)^2}-1\Big)=0\,.$$ 
In particular, there exists $i_0$ such that 
$\lb'_j(e'_j\cdot e_{i_0}^*)^2\ge\lb'_k(e'_k\cdot e_{i_0}^*)^2$. 
\end{lemma} 

\begin{proof} 
We shall prove a slightly more general result. 
By Lemma~\ref{lemAMB2}, (3), we have 
$\sum_i\abs{e'_k\cdot e_i^*}\,\abs{e_i\cdot{e'_k}^*}=1$ 
and 
$\sum_i\abs{e'_j\cdot e_i^*}\,\abs{e_i\cdot{e'_j}^*}=1$. 
Using Lemma~\ref{lembases}, (4), we obtain for all $i$ such that 
$e'_k\cdot e_i^*\ne 0$, 
$$\abs{e'_j\cdot e_i^*}\,\abs{e_i\cdot{e'_j}^*}= 
\abs{e'_k\cdot e_i^*}\,\abs{e_i\cdot{e'_k}^*}\, 
\frac{\lb'_j}{\lb'_k}\,
\frac{(e'_j\cdot e_i^*)^2}{(e'_k\cdot e_i^*)^2}\,.$$ 
This implies the identity 
$$\sum_{i,\,e'_k\cdot e_i^*\ne 0} 
\abs{e'_k\cdot e_i^*}\,\abs{e_i\cdot{e'_k}^*}\, 
\Big(\frac{\lb'_j}{\lb'_k}\,
\frac{(e'_j\cdot e_i^*)^2}{(e'_k\cdot e_i^*)^2}-1\Big)+ 
\sum_{i\mid e'_k\cdot e_i^*=0}
\abs{e'_j\cdot e_i^*}\,\abs{e_i\cdot{e'_j}^*}=0\,,$$ 
which completes the proof of the lemma. 
\end{proof} 

In the next lemma, we exceptionnally forget for the Euclidean 
structure considering a more general bilinear form, 
still denoted by $x\cdot y$. 

\begin{lemma}\label{lemformal} 
Let $K$ be a field of characteristic not~$2$ and let $V$ 
be an $n$-dimensional vector space over~$K$, equipped with a basis 
$\cB_0=(\va_1,\dots,\va_n)$. 
Denote by $R$ the ring $K[t_{i,j}]$, $1\le i\le j\le n$ 
of polynomials in $\frac{n(n+1)}2$ variables $t_{i,j}$. 
Consider on $V$ the symmetric bilinear form with values in~$R$ 
such that $e_j\cdot e_i=e_i\cdot e_j=t_{i,j}$ for $1\le i\le j\le n$. 
Let $(e_1,\dots,e_p)$ and $(e'_1,\dots,e'_{p'})$ be two systems 
of vectors of $V$. Let $\lb_i,\,\lb'_j$, $1\le i\le p$, 
$1\le j\le p'$ be elements of $K$. Then the relations 
$$\sum_{i=1}^p\,\lb_i\,N(e_i)= \sum_{j=1}^{p'}\,\lb'_j\,N(e'_j) 
\ \nd\ 
\sum_{i=1}^p\,\lb_i\,N(e_i)\,p_{e_i}= 
\sum_{j=1}^{p'}\,\lb'_j\,N(e'_j)\,p_{e'_j}$$ 
are equivalent. 
{\small\rm 
[As above, $N(x)=x\cdot x$, and $p_x$ denotes the orthogonal projection 
to~$x$ with respect to the given bilinear form.]} 
\end{lemma} 

\begin{proof} 
For convenience, we set $t_{j,i}=t_{i,j}$ if $j<i$. 
Write $e_i=\sum_k\,x_{k,i}\va_k$ and $e'_j=\sum_k\,x'_{k,j}\va_k$. 
The equation 
$\sum_{i=1}^p\,\lb_i\,N(e_i)= \sum_{j=1}^{p'}\,\lb'_j\,N(e'_j)$ 
reads 
$\sum_{i,k,\ell}\,\lb_i\,x_{\ell,i}\,x_{k,i}\,t_{\ell,k}= 
\sum_{j,k,\ell}\,\lb'_j\,x'_{\ell,j}\,x'_{k,j}\,t_{\ell,k}$, 
equivalent to 
$$\forall\,k,\ell,\ \, 
\sum_i\,\lb_i\,x_{\ell,i}\,x_{k,i}= 
\sum_j\,\lb'_j\,x'_{\ell,j}\,x'_{k,j}\,.$$ 

Identifying the matrices with respect to $\cB_0$ of both sides 
of the first equality, we obtain 
$\forall\,k,m,\,
\sum_i\,\lb_i\,x_{k,i}\,x_{\ell,i}t_{\ell,m}= 
\sum_j\,\lb'_j\,x'_{k,j}\,x'_{\ell,j}\,t_{\ell,m}\,.$ 
Equating the coefficients of the variables $t_{\ell,m}$, 
we recover the equality displayed above. 
\end{proof}

\section{Perfection relations for lattices.}\label{seclatperfrel} 

We now consider perfection relations of the form 
$$\sum_{i=1}^n\,\lb_i\,p_{e_i}=\sum_{i=1}^n\,\lb'_i\,p_{e'_i} 
\eqno{(\dag)}$$ 
where the $e_i$ and the $e'_i$ constitute two rank~$n$ sets 
of minimal vectors in a lattice and the coefficients $\lb_i,\lb'_i$ 
are strictly positive. We denote by $\Lb_0$, $\Lb'_0$, $\Lb$ 
the lattices generated by the vectors $e_i$, the vectors $e'_j$, 
and all vectors $e_i,e'_j$ respectively. 

The quotient $\Lb/\Lb_0$ defines a code $C$ over $\Z/d\Z$, 
where $d$ denotes the annihilator of $\Lb/\Lb_0$. 
We define similarly $d'$ and $C'$ with respect to $\Lb'_0$. 

We say that the perfection relation is \emph{regular} if $d'=d$ 
and the two codes are equivalent, and \emph{irregular} otherwise. 

\begin{lemma}\label{lemvmin} 
Consider in some lattice $L$ a relation 
$\sum_{i=1}^k\,\lb_i\,N(e_i)=\sum_{i=1}^{k'}\,\lb'_i\,N(e'_i)$ 
with real coefficients $\lb_i,\,\lb'_j$ such that 
$\sum_i\,\lb_i=\sum_j\,\lb'_j$ and non-zero vectors $e'_j$. 
If the $e_i$ are minimal and the $\lb'_j$ are strictly positive, 
then the $e'_j$ are also minimal. 
\end{lemma} 

\begin{proof} 
Set $m=\min L$. We have 
$$0=-\sum_i\,\lb_i\,m+\sum_j\,\lb'_j\,N(e'_i)= 
\sum_j\,\lb'_j\big(N(e'_j)-m\big)\,.$$ 
Since all terms in the second sum are non-negative, 
all must be zero. 
\end{proof} 

We say that the lattice $\Lb$ is \emph{\pfi} 
if the system $\{e_i,e'_j\}$ is irreducible in the sense 
of Definition~\ref{defirred}. 
A lattice endowed with a perfection relation as above 
is in a unique way a direct sum (not necessarily orthogonal) 
of {\pfi} sublattices. 

\smallskip 

Our two main sources of perfection relations are: 

\noi(1) Relations of the form 
$\sum_{i=1}^n\,p_{e_i}=\sum_{i=1}^n\,p'_{e_i}$ when $(e_i)$ and $(e'_i)$ 
are orthogonal bases for $E$; the relation above holds because both sides 
are equal to the identity; 

\noi(2) Relations which stem directly or not from ``Watson's condition'', 
that we define below. 

\noi[However, other kinds of perfection relations exist 
in dimensions $n\ge 8$, see below.] 

\smallskip 

We first state an identity due to Watson whose proof is left to the reader. 

\begin{lemma}\label{lemwatid}\emph{(Watson.)} 
Let $f_1,\dots,f_n$ be $n$~independent vectors in $E$, 
and let $f=\dfrac{\sum_{i=1}^\ell\,a_i\,f_i}d$
where $\ell\ge 2$ is an integer and $d>0$ and $a_i$ are real numbers. 
Denote by $\sgn(x)$ the sign of the real number $x$ 
($\sgn(x)=0$ if $x=0$). Then 
$$\big((\sum_{i=1}^\ell\abs{a_i})-2d\big)N(f)= 
\sum_{i=1}^\ell\abs{a_i}\big(N(f-\sgn(a_i)f_i)-N(f_i)\big)\,. 
\eqno{(\ast)}\qed $$ 
\end{lemma} 

The first part of the following lemma is due to Watson: 

\begin{lemma}\label{lemwat} 
Assume now that the $f_i$ are independent minimal vectors 
in a lattice $L$ containing also $f$ 
and that $d>1$ and $a_1,\dots,a_\ell$ are non-zero integers. 
Set $A=\sum_{i=1}^\ell\abs{a_i}$. Then we have $A\ge 2d$, 
and equality holds \ifff all vectors $f-\sgn(a_i)\,f_i$ are minimal. 

Moreover, when equality holds 
(we then say that \emph{Watson's condition} holds), 
we have $\abs{a_i}\le\frac d2$ for all $i$, 
and if $d\ge 4$, equality holds for at most one index~$i$. 
\end{lemma} 

\begin{proof} 
The vectors $f-\sgn(a_i)\,f_i$ are non-zero vectors in $L$ 
(because $d>1$). Hence the Right Hand Side of $(\ast)$ 
is non-negative, and is zero \ifff all vectors $f-\sgn(a_i)\,f_i$ 
are minimal. 

Suppose now that $A=2d$. Since we also have 
$\sum_{i\ne i_0}\abs{a_i}+\abs{d-\abs{a_{i_0}}}\ge 2d$ 
for all $i_0$ 
(consider $f-\sgn(a_{i_0})\,f_{i_0}$ instead of $f$), 
we have $d-2\abs{a_{i_0}}\ge 0$. 

For the remaining of the proof, we assume that all $a_i$ are positive 
(we easily reduce to this case by negating some $f_i$ if necessary). 
Let $k$ be the number of $a_i$ with $a_i=\frac d2$. 
If $k>0$, $d$ is even, say, $d=2d'$, and 
$$2f-\sum_{a_i=d'}f_i= 
\frac{\sum_{a_i<d'}a_i\,f_i}{d'}$$ 
satisfies the first part of the lemma, 
i.e. $\sum_{a_i<d'}\,a_i\ge 2d'$, hence 
\linebreak 
$2d'+k d'\le 4d'$, i.e. $k\le 2$. 

Finally, assume that $k=2$ and, say, that we have 
$a_{\ell-1}=a_\ell=d'$ and $a_i<d'$ for $i\le \ell-2$. 
Let $f'=\dfrac{\sum_{i=1}^{\ell-2}\,a_i\,f_i}{d'}$. 
We have $\sum_{i=1}^{\ell-2}\,a_i=2d'$, 
hence the vectors $f'-f_i$ are minimal for $i\le\ell-2$ 
(because $d'>1$), 
and since $f=\frac{(f'-f_i)+f_i+f_{\ell-1}+f_\ell}2$, 
the four vectors of the numerator are mutually orthogonal, 
which shows that $N(f')=2$ and $f_i\cdot f'=1$, 
thus that $N(f-f_i-f_j)=2\,f_i\cdot f_j$ for $1\le i<j\le\ell-2$. 
Thus we have $f_i\cdot f_j\ge\frac 12$, 
hence $f_i\cdot f_j=\frac 12$, and 
$$f_i\cdot f'=\frac{\sum_{j=1}^{\ell-2} a_j\,(f_i\cdot f_j)}{2d'}= 
\frac{a_i}{2d'}+\frac{\sum_{j=1}^{\ell-2} a_j}{2d'}= 
1+\frac{a_i}{2d'}\,,$$ 
a contradiction. 
\end{proof} 

\begin{lemma}\label{lemsumxi} 
Let $L$ be a lattice, let $f_1,\dots,f_n$ be independent minimal 
vectors of $L$ and let $x=\sum_i\,x_i\,f_i\in L$, $x$ non-zero 
and different from the $\pm f_i$. 
Then $\sum_i\,\abs{x_i}\ge 2\,.$ 
\end{lemma} 

\begin{proof} 
If $x$ does not belong to the lattice generated by $f_1,\dots,f_n$, 
this follows from Lemma~\ref{lemwat}. 
Otherwise, the $x_i$ are integers, and at least two are non-zero. 
\end{proof} 

The proposition below is an immediate consequence 
of Lemma~\ref{lemformal}; for ther sake of simplicity, 
we assume that the $a_i$ are strictly positive. 

\begin{prop}\label{propwat} 
Assume that Watson's condition holds in Lemma \ref{lemwat}, 
Then the $2\ell$ vectors $f_i,f-a_i\,f_i$ 
satisfy the (unique up to proportionality) perfection relation 
$$\sum_{i=1}^\ell\,a_i\,p_{f_i}= 
\sum_{i=1}^\ell\,a_i\,p_{f-a_i\,f_i}\,. 
\eqno{(\ast\,\ast)}\qed$$ 
\end{prop} 

We now return to the notation $\Lb,\Lb_0,\Lb'_0$ 
introduced at the beginning of this section. 
The following lemma, that we prove using Watson's inequality, 
gives a necessary condition for the system $(\Lb,\Lb_0,\Lb'_0)$ 
to exist. 

\begin{lemma}\label{lemlength} 
Let $f_0,f_1,\dots,f_r$ be a system of representatives of $\Lb$ 
modulo $\Lb_0$ with $f_0=0$. For every~$i\ge 1$, 
write $f_i=\sum_j\,x_j^i\,e_j$. Then for all~$j$, 
there exists $i\ge 1$ such that $x_j^i\notin\Z$. 
\end{lemma} 

\begin{proof} 
Suppose there exists some $j$ such that $\forall\,i,\,x_j^i\in\Z$. 
Then, every $x\in\Lb$ has an integral component 
$x\cdot e_j^*$ on~$e_j$. Let 
$$K=\{k\mid e'_k\cdot e_j^*\ne 0\}\,.$$ 
For $k\in K$, we have $\abs{e'_k\cdot e_j^*}\ge 1$, 
hence Lemma~\ref{lemAMB2},\,(2) implies the upper bound 
$\sum_{k\in K}\,\abs{e_j\cdot{e'_k}^*}\le 1$. 
By Lemma~\ref{lemAMB2},\,(1), we also have 
$\sum_{k=1}^n\,\abs{e_j\cdot{e'_k}^*}\le 1$, 
which contradicts Lemma~\ref{lemsumxi}. 
\end{proof} 

\begin{corol}\label{corindex} 
The lattice $\Lb$ contains strictly $\Lb_0$ and $\Lb'_0$. 
\end{corol} 

\begin{proof} 
The strict inclusion $\Lb\supset_{\ne}\Lb_0$ follows immediately 
from the lemma above. Exchanging $\Lb_0$ and $\Lb'_0$ 
proves the second one. 
\end{proof} 

\begin{corol}\label{corcyclic} 
If $\Lb/\Lb_0$ is cyclic, say, $\Lb=\la\Lb_0,e\ra$ 
for $e=\frac{a_1 e_1+\dots+a_n e_n}d$ with $0\le\abs{a_i}\le \frac d2$,
then all $a_i$ are non-zero. 
\end{corol} 

\begin{proof} 
The proof is immediate. 
\end{proof} 

\begin{notation} 
With the notation of Lemma~\ref{lemwat}, we may assume 
that $\abs{a_i}\le\frac d2$. For $i=1,2,\dots,d'=\lf\frac d2\rf$, 
we then denote by $m_i$ the number of coefficients $a_j$ equal 
to $\pm i$. 
\end{notation} 

With this notation, we have $m_1+\dots+m_{d'}=\ell$ 
and Watson's inequality reads $m_1+2m_2+\dots+d' m_{d'}\ge 2d$. 

\smallskip 

Other identities involving systems of $2n$~vectors of rank~$n$ 
exist for~$n\ge 8$. 

We first consider the case when $d=5$, setting 
$f'=2f-\sum_{i=m_1+1}^\ell\,f_i$, $f'_i=f'-f_i$ if $i\le m_1$ 
and $f'_i=f'-f_i$ if $i>m_1$. Then we have the formal identity 
$$\sum_{i=1}^\ell\big(N(f'_i)-N(f_i)\big)=(m_2-4)\,N(f)+(m_1-4)\,N(f')\,,$$ 
an identity considered by Zahareva when $(m_1,m_2)=(4,3)$ or $(4,4)$. 
In this last case, Lemmas \ref{lemformal} and \ref{lemvmin} show 
that when the $f_i$ are minimal, the $f'_i$ 
are also minimal and that $\sum_i\,p_{f_i}=\sum_i\,p_{f'_i}$. 

An identity of the same kind exists for $d=7$, involving 
the reductions modulo $\Lb_0$ of $2f$ and $3f$ and the differences 
$m_i-3$, and when $m_1=m_2=m_3=3$ (hence $\ell=9$), 
there again exists a perfection relation as above. 
An example is given in \cite{M1}, Rem.~9.2. 

\smallskip 

In the following lemma, we compare the orders of the various basis 
vectors $e'_j$ modulo $\Lb_0$. 

\begin{lemma}\label{lemAMB4} 
For $d>0$, let $N_d$ be the set of vectors $e'_j$ of order $d$ 
in the quotient $\Lb/\Lb_0$ and let $\nu_d$ be the number 
of such vectors. Then we have $\nu_1\le\sum_{d\ge 3}\,(d-2)\nu_d$, 
and when equality holds, we have $\abs{e'_j\cdot e_i^*}=\frac 1d$ 
for all $e'_j$ of order~$d$ and all $i$ 
such that $e'_j\cdot e_i^*\ne 0$. 
\end{lemma} 

\begin{proof} 
First consider a vector $e'_j\in N_d$. 
For every $i$ such that $e_i\cdot{e'_j}^*\ne 0$, 
$e'_j\cdot e_i^*$ is non-zero by Lemma~\ref{lemAMB2},\,(1), 
hence $\ge\frac 1d$. 
By Assertion (3) of the same lemma, we have 
$$1=\sum_{i,\,e_i\cdot {e'_j}^*\ne 0}\, 
\abs{e_i\cdot {e'_j}^*}\abs{\,e'_j\cdot e_i^*}\ge\frac 1d\, 
\sum_{i,\,e_i\cdot {e'_j}^*\ne 0}\, \abs{e_i\cdot {e'_j}^*}\,,$$ 
hence $\sum_{i=1}^n\,\abs{e_i\cdot {e'_j}^*}\le d$. 
Taking now the sum over $j$, we obtain 
$\sum_{1\le i,j\le n}\,\abs{e_i\cdot {e'_j}^*}\le\sum_d\,d\nu_d$. 
By Lemma~\ref{lemsumxi}, we have 
$\sum_j\,\abs{e_i\cdot {e'_j}^*}\ge 2$ for all $i$, 
which implies $2n\le\sum_d\,d\nu_d$, i.e. $\sum_d\,(d-2)\nu_d\ge 0$ 
(because $n=\sum_d\nu_d$). This completes the proof of the lemma. 
\end{proof}

\section{2-elementary quotients.}\label{seclat2elem} 

We keep the notation of Section~\ref{seclatperfrel}, 
and consider perfection relations on lattices $\Lb$ 
such that $\Lb/\Lb_0$ is $2$-elementary. 
We first construct some examples concerning root lattices. 

\smallskip 

Recall that $\D_n$ is the even sublattice of the lattice $\Z^n$, 
endowed with its canonical basis $(\va_1,\dots,\va_n)$. 
It is generated by its roots $\pm(\va_i\pm\va_j)$. 
One easily checks that orthogonal frames (of minimal vectors) 
exist \ifff $n$ is even and $n\ge 4$, and that they then constitute 
a unique orbit under the action of the Weyl group $W(\D_n)$, 
namely that of 
$e_1=\va_1+\va_2$, $e_2=\va_1-\va_2$, $e_3=\va_3+\va_4$, \dots, 
$e_n=\va_n-\va_{n-1}$. For any automorphism $\s$ of $\Z^n$, 
we have the equalities $\sum_i\,p_{e_i}=\Id=\sum_i\,p_{\s(e_i)}$. 
Choosing $\s$ such that the lines $\R\,\s(e_i)$ are distinct 
from the lines $\R\,e_j$, we obtain a perfection relation. 
A possible choice for $\s$ is the product of transpositions 
$(2,3)(4,5)\dots(n-2,n-1)$. 

\smallskip 

The same kind of result works for $\E_8 =\D_8\cup(e+\D_8)$ 
with $e=\frac{\va_`+\dots+\va_8}2$ and for $\E_7$ 
(the orthogonal complement of a root in~$\E_8$). 
The group $\Aut(\E_8)=W(\E_8)$ acts transitively on the pairs 
of orthogonal roots, which can be taken inside $\D_8$, 
whose orthogonal complement is then isometric to $\D_6$. 
Hence there is again a unique orbit of orthogonal frames 
of minimal vectors in $\E_8$, and then also in $\E_7$. 

\smallskip 

Finally, 
\emph{there exist in the lattices $\D_n$ ($n\ge 4$ even), 
$\E_7$ and $\E_8$ perfection relations of the form 
$\sum\,p_{e_i}=\sum\,p_{e'_i}$ with mutually orthogonal 
systems $(e_i)$ and $(e'_i)$. In particular, $\Lb/\Lb_0$ 
and $\Lb/\Lb'_0$ are isomorphic, $2$-elementary groups, of order 
$2^{(n-2)/2}$, $2^3$ and $2^4$ respectively.} 

\smallskip 

Besides these regular relations, there exists an irregular one 
in the case of $\E_8$, that we now construct. 

Start with an orthogonal frame $(e_1,\dots,e_8)$ of roots. 
For $i=1,3,5,7$, consider the mutually orthogonal, norm~$6$ vectors 
{\small 
$$f_1=e_4+e_6+e_8,\,f_3=e_2-e_4+e_8,\,f_5=e_2+e_4-e_6,\,f_7=e_2+e_6-e_8$$} 
(where the signs have been chosen in accordance with the construction 
of the ternary tetracode). Define the vectors $e'_i$ by 
{\small 
$$e'_1=\frac{e_1+f_1}2,\ e'_2=\frac{-e_1+f_1}2,\ e'_3=\frac{e_3+f_3}2,\dots, 
\ e'_8=\frac{-e_7+f_7}2\,.$$} 
For $i=1,3,5,7$, we have $e'_i\cdot e'_i=2$ and $e'_i\cdot e'_{i+1}=1$, 
so that $e'_i$ and $e'_{i+1}$ generate a hexagonal lattice, 
whose minimal vectors are (up to sign) $e'_i$, $e'_{i+1}$ and $e_i$. 
The planes $H_i$ of these lattices are mutually orthogonal, 
so that $\Id=p_{H_1}+p_{H_3}+p_{H_5}+p_{H_7}$, 
and for every~$i$, we have 
$p_{H_i}=\frac 23\,(p_{e_i}+p_{e'_i}+p_{e'_{i+1}})$. 
This implies the irregular relation 
$$\sum_{i\text{ odd}}\,p_{e_i}+3\,\sum_{i\text{ even}}\,p_{e_i}= 
2\,\sum_{j=1}^8\,p_{e'_j}\,.$$ 
Here, $\Lb/\Lb_0$ is $2$-elementary of order $16$ whereas 
$\Lb/\Lb'_0$ is $3$-elementary of order $9$, as one sees 
writing 
{\small 
$$e_2=\frac{e'_3+e'_4+e'_5+e'_6+e'_7+e'_8}3\ \nd\ 
e_4=\frac{e'_1+e'_2-e'_3-e'_4+e'_5+e'_6}3\,.$$} 

\begin{theorem}\label{th2elem} 
Let $\Lb$ be a {\pfi} lattice 
which possesses a perfection relation 
$\sum_i\,\lb_i\,p_{e_i}=\sum_i\,\lb'_i\,p_{e'_i}$ with strictly 
positive coefficients $\lb_i,\,\lb'_i$ and independent 
systems $(e_i)$ and $(e'_i)$ of minimal vectors. 
Assume moreover that $\Lb/\Lb_0$ is $2$-elementary. 
Then $\Lb$ is similar to one of the root lattices 
$\D_n$ ($n\ge 4$ even), $\E_7$, $\E_8$ endowed with a regular 
relation, or to $\E_8$, endowed with the irregular one. 
\end{theorem} 

The proof of this theorem will occupy the remaining of the section. 
We consider lattices $\Lb$, $\Lb_0$ as above, scaled to minimum~$1$ 
as in Section~\ref{seceuclid}, and assume from now on 
that $\Lb/\Lb_0$ is $2$-elementary. 
We begin with three lemmas. 
The notation $A_i$, $A'_i$ is that of Lemma~\ref{lemAMB1}. 

\begin{lemma}\label{lem2elem1} 
All $e'_i$ are of the form 
$e'_i=\dfrac{\pm e_{j_1}\pm e_{j_2}\pm e_{j_3}\pm e_{j_4}}2$, 
where the vectors $e_{j_1},e_{j_2},e_{j_3},e_{j_4}$ 
are mutually orthogonal. 
\end{lemma} 

\begin{proof} 
If $e'_i$ belongs to $\Lb_0$, we have $A'_i\le -1$ 
(see Lemma~\ref{lemsumxi}). Otherwise, $e'_i$ is of the form 
$\frac{a_1\,e_1+\dots+a_n e_n}2$ with $k\ge 4$ 
odd coefficients $a_i$. 
Then we have $A'_i=1-\frac{\sum_i\,a_i^2}4\le 1-\frac k4\le 0$. 
By Lemma~\ref{lemAMB1}, all $A'_i$ must be zero, 
which is possible only if no $e'_i$ lies in $\Lb_0$, 
and if for all $e'_i$, $k=4$ and all $a_i\ne 0$ are $\pm 1$. 
Finally, we recognize the configuration of a $4$-dimensional 
centred cubic lattice. 
\end{proof} 

\begin{lemma}\label{lem2elem2} 
\begin{enumerate} 
\item 
For all $j$, 
$\sum_{i\in\supp{e_j}}\,\abs{e_j\cdot {e'}_i^*}=2$. 
\item 
For all $i$, 
$\sum_{j\in\supp{e'_i}}\,\abs{e_j\cdot {e'}_i^*}=2$. 
\end{enumerate} 
\noi
Moreover, in {\rm(2)}, either the four terms $\abs{e_j\cdot {e'}_i^*}$ 
are equal to $\frac 12$, or exactly one of them is equal to~$1$, 
and then the corresponding $e_j$ belongs to~$\Lb'_0$. 
\end{lemma} 

\begin{proof} 
\noi{\sl Proof.} 
(1) and (2) result from Lemma~\ref{lemAMB2}, observing that 
$$i\in\supp e_j\Equi j\in\supp e'_i\Equi 
e'_i\cdot e_j^*=\pm\frac12\,.$$ 
Moreover, (1) tells us that Watson's condition is satisfied 
by any $e_j$ not in $\Lb'_0$, and Lemma~\ref{lemwat} 
shows that the components of such an $e_j$ are bounded above 
by $\frac 12$. This completes the proof of the lemma. 
\end{proof} 

\begin{lemma}\label{lem2elem3} 
If no vector $e_i$ belongs to $\Lb'_0$, $\Lb/\Lb'_0$ 
is also $2$-elementary, every $e_i$ is of the form 
$e_i=\dfrac{\pm e'_{j_1}\pm e'_{j_2}\pm e'_{j_3}\pm e'_{j_4}}2\,,$ 
and two vectors $e'_{i_1}$, $e'_{i_2}$ with the same support 
differ by an even number of minus signs. 
\end{lemma} 

\begin{proof} 
For any pair $(i,j)$, by Lemma~\ref{lem2elem2}, 
we have $e_i\cdot{e'}_j^*\in\{0,\pm\frac 12\}$, 
which implies that $A_j=1-\sum_i\,(e_j\cdot{e'}_i^*)^2\le 0$. 
By Lemma~\ref{lemAMB1}, all $A_j$ must be zero, 
which completes the first part of the proof. 

If two vectors $e'_{i_1}$, $e'_{i_2}$ with the same support 
differ by an odd number of minus signs, 
then $e'_{i_1}+e'_{i_2}$ or $e'_{i_1}-e'_{i_2}$ is of the form 
$\pm e_j$, a contradiction. 
\end{proof} 


\noi{\it Proof of Theorem~\ref{th2elem} in the case when no vector 
$e_i$ belongs to~$\Lb'_0$.} 

Observe that from Lemmas \ref{lemAMB2} and \ref{lem2elem3}, 
there exist for all $j$ exactly four indices $i$ such that 
$e'_i\cdot e_j^*\ne 0$, and that $e_j$ is of the form 
$$e_j=\frac{\sum_i\sgn(e'_i\cdot e_j^*)\,e'_i}2\,.$$ 

Suppose first that two vectors $e'_i$ have the same support, 
say (thanks to Lemma~\ref{lem2elem3}), 
$e'_1=\frac{e_1+e_2+e_3+e_4}2$ and $e'_2=\frac{-e_1-e_2+e_3+e_4}2$. 
There exist two more vectors $e'_3,e'_4$ having $e_3$ 
in the numerator with, say, $e'_3\cdot e_3^*$ and $e'_3\cdot e_4^*$ 
positive. The displayed formula above shows 
that $e'_3+e'_4=e_3-e_4$. Then $e'_3$ and $e'_4$ also have 
the same support $\{3,4,i,j\}$. Since two supports cannot have 
three common indices (otherwise, the code defined by the numerators 
of the $e'_i$ would have weight~$2$), 
the pair $\{i,j\}$ is either equal to 
or disjoint from $\{1,2\}$. In the first case, $\Lb$ contains 
the centred cubic lattice as a component, 
hence is similar to $\D_4$ since it is {\pfi}. 
In the second case, write $(i,j)=(5,6)$. Applying the preceding 
argument to $\{5,6\}$, we construct a sequence 
$\{5,6\}$, $\{7,8\},\,\dots$, which must break down 
when the last $\{2p-1,2p\}$ meets $\{1,2\}$, 
since $e_1$ and $e_2$ must occur in four vectors~$e'_j$. 
Since $\Lb$ is {\pfi}, $\dim\Lb=2p$. 

Any sum of distinct vectors $e'_i$ is congruent modulo $\Lb_0$ 
to a vector of the form $e=\frac{e_{i_1}+e_{i_2}+e_{i_3}+e_{i_4}}2$, 
and any pair $(e_i,e_j)$ occurs in the numerator of such a vector. 
Hence all $e_i$ are mutually orthogonal, 
since $e,e_{i_1},\dots,e_{i_4}$ generate a centred cubic lattice. 
Similarly, the $e'_j$ are mutually orthogonal. 
Since all scalar products $e_i\cdot e'_j$ are equal 
to $0$ or $\pm\frac12$, the lattice $\Lb$ rescaled to minimum $2$ 
is an irreducible root lattice, which contains $\Lb_0$ 
(of determinant $2^n$) to index $2^{(n-2)/2}$. 
Therefore $\det(\Lb)=4$, hence $\Lb$ is similar to $\D_{2p}$, 
endowed with a perfection relation coming from orthogonal systems 
of minimal vectors, i.e. a regular one. 

\smallskip

Next suppose that two distinct vectors $e'_i$ have distinct supports. 
A ``sudoku-like'' proof will allow us to conclude. 

Using the displayed formula above, we show that we may take 
the first four vectors $e'_1,\dots,e'_4$ of the form 
{\small 
$$\frac{e_1+e_2+e_3+e_4}2,\, 
\frac{e_1-e_2+e_5+e_6}2,\, 
\frac{e_1-e_3-e_5+e_7}2,\, 
\frac{e_1-e_4-e_6-e_7}2\,.$$ 
} 
Then $e_2$ shows up in, say, $e'_5$ and $e'_6$, 
say, $e'_5=\frac{e_2-e_3+e_j\pm e_k}2$ with 
\linebreak 
$(j,i)=(6,7)$ or $(5,8)$. 

In the first case, we check that the last three vectors are 
{\small 
$$\frac{e_2-e_3+e_6-e_7}2,\, 
\frac{e_2-e_4+e_5+e_7}2,\, 
\frac{e_3-e_4-e_5+e_6}2\,.$$} 

In the second case, we check that the last four vectors are 
{\small 
$$\frac{e_2-e_3+e_5+e_8}2,\, 
\frac{e_2-e_4+e_6-e_8}2,\, 
\frac{e_3-e_4+e_7+e_8}2,\, 
\frac{e_5-e_6+e_7-e_8}2\,.$$} 

In both cases, checking the scalar products and calculating 
the determinant as we did for $\D_n$ shows that the lattice $\Lb$ 
is similar to $\E_7$ and $\E_8$ respectively, 
endowed with the regular perfection relation. 

\medskip 

\noi{\it Proof of Theorem~\ref{th2elem} in the case when at least %
one vector $e_i$ belongs to~$\Lb'_0$.} 

By Lemma~\ref{lem2elem2},\,(2), such a vector $e_i\in\Lb'_0$ 
is of the form $e_i=\pm e'_j\pm e'_k$, and if another $e_{i'}$ 
is also of the form $\pm e'_{j'}\pm e'_{k'}$, then $j',k'$ 
are distinct from $j,k$. Then up to permutation and sign changes, 
we may write $e_i=e'_i-e'_{i+1}$ with $i$ odd. 

Let $j\ne i$ with $i\in\supp e_j$, and write 
$e_j=\sum_{k=1}^n\,x_k\,e'_k$. We have $x_{i+1}=x_i$, 
and in particular, $i+1\in\supp e_j$: 
we have $e'_k\cdot e_i^*=0$ or $\pm\frac 12$, 
and its vanishing and sign is given by Lemma~\ref{lemAMB2}, 
hence $0=e_j\cdot e_i^*=\frac{x_i}2-\frac{x_{i+1}}2\,.$ 
The vector $e_j$ is not in $\Lb'_0$ (because $j\ne i$), 
and cannot have $2$ as a denominator (since it satisfies Watson's 
condition, it would be a sum of $4$~mutually orthogonal vectors, 
among which $e'_i$ and $e'_{i+1}$, whereas 
$e'_i\cdot{e'}_{i+1}=\frac 12\ne 0$). Hence $e_j$ has denominator 
$d\ge 3$. We now show that at least six of its components 
have absolute values strictly smaller than $\frac 12$: 
if its support has $\ell\ge 8$ elements or if $d$ is odd, 
this comes from Lemma~\ref{lemwat}. Otherwise, 
we have $d\ge 4$ and $\ell=7$ (see \cite{M1}), in which case 
six out of the seven components are equal to $\frac 14$. 

\smallskip 

We now consider the (at least) six ``small'' non-zero 
components of~$e_j$. Let $x_k$ such a component. 
By Lemma~\ref{lem2elem2}, there exists $j'\in\supp e'_k$ 
such that $e_{j'}$ belongs to $\Lb'_0$, actually in our notation, 
$j'=k$ or $k-1$, say, $j'=k$, hence $e_k=e'_k-e'_{k+1}$. 
We may thus write 
$$e_j=x_i(e'_i+e'_{i+1})+x_k(e'_k+e'_{k+1})+ 
x_\ell(e'_\ell+e'_{\ell+1})+\dots\,,$$ 
with $x_i,x_k,x_\ell\ne 0$. 

\smallskip 

We now go into the proof, taking precise notation as follows: 
{\small 
$$e_1=e'_1-e'_2\,,\text{ with } 
e'_1=\frac{e_1+e_2+e_4+e_6}2\ \nd\ e'_2=\frac{-e_1+e_2+e_4+e_6}2\,,$$} 
and writing $e_2,e_4,e_6$ as above: 
{\small 
$$e_2=x_1(e'_1+e'_2)+x_3(e'_3+e'_4)+x_5(e'_5+e'_6)+ \dots$$} 
with $x_1,x_3,x_5\ne 0$, and similarly 
{\small 
$$e_4=y_1(e'_1+e'_2)+y_3(e'_3+e'_4)+y_5(e'_5+e'_6)+ \dots\,$$} 
and 
{\small 
$$e_6=z_1(e'_1+e'_2)+z_3(e'_3+e'_4)+z_5(e'_5+e'_6)+ \dots\,.$$} 
Here, $x_1,y_1,z_1$ are strictly positive by Lemma~\ref{lemAMB2}, 
and up to a global change of signs in the relations 
$e_3=e'_3-e'_4$, $e_5=e'_5-e'_6$, 
we may assume that $x_3$ and $x_5$ are positive. 
Using the relation $e'_1+e'_2=e_2+e_4+e_6$, we obtain 
$$x_1+y_1+z_1=1,\,x_3+y_3+z_3=0,\,\nd\,x_5+y_5+z_5=0\,.\eqno{(\ast)}$$ 
Since the support of $e'_3$ is distinct from that of $e'_1$ 
(because $3$ belongs to it), it does not contain $\{2,4,6\}$, 
which implies $x_3y_3z_3=x_5y_5z_5=0$. 
We may clearly suppose that $z_3=0$ (thus $y_3=-x_3<0$ by $(\ast)$). 
Let us prove that $y_5=0$. We may write 
$e'_3=\frac{e_3+e_2-e_4+e_8}2$ for some $e_8$ 
(see again Lemma~\ref{lemAMB2}) whose support does not contain 
$1$ and~$2$, say 
$$e_8=t_3(e'_3+e'_4)+t_5(e'_5+e'_6)+
t_7(e'_7+e'_8)+ 
\dots\,,$$ 
with $t_3>0$ by the same lemma, thus $x_3y_3t_3\ne 0$, 
which implies as above $x_5y_5t_5=0$, and $t_7>0$ for convenience. 
Using the relation $e'_3+e'_4=e_2-e_4+e_8$, we obtain 
$$x_1-y_1=0,\,x_3-y_3+t_3=1,\,\nd\,x_5-y_5+t_5=0\,.\eqno{(\ast\ast)}$$ 
If $y_5$ were non-zero, we would have $z_5=0$, then $y_5=-x_5$ 
by $(\ast)$ and $t_5=-2x_5$ by $(\ast\ast)$ would be non-zero, 
a contradiction. We thus have $y_5=0$ and $z_5=t_5=-x_5<0$, 
hence $e'_5=\frac{e_5+e_2-e_6-e_8}2$. Denote by $x_7,y_7,z_7$ 
the components of $e_2,e_4,e_6$ respectively on $e'_7+e'_8$. 
We have as above $x_7y_7t_7=x_7z_7t_7=0$ and also $x_7+y_7+z_7=0$ 
(using again the relation $e_2+e_4+e_6=e'_1+e'_2$). 
Since $t_7>0$, $x_7$ must vanish. 

Using all components of $e'_i+e'_{i+1}$, $i=1,3,5,7$, 
on $e_2,e_4,e_6,e_8$, we obtain 
{\small 
$$x_1=y_1=z_1=\,x_3=-y_3=t_3=\,x_5=-z_5=-t_5=\,y_7=-z_7=t_7= 
\frac 13\,.$$} 
We know for $e_2,e_4,e_6,e_8$ $6$~components equal to $\pm\frac 13$. 
Since they satisfy Watson's condition, they do not have 
any other non-zero component, hence are uniquely defined, 
and since $\Lb$ is assumed to be \pfi, it has dimension~$8$, 
and we recognize the irregular perfection relation for $\E_8$. 
This completes the proof of Theorem~\ref{th2elem}. \qed

\section{Lattices of index~3}\label{secindex3} 

We keep the notation of Section~\ref{seclatperfrel}, 
and now study perfection relations when $[\Lb:\Lb_0]=3$. 

\begin{theorem}\label{thindex3} 
Let $\Lb$ be a lattice endowed with a perfection relation 
$\sum_{i=1}^n\,\lb_i\,p_{e_i}=\sum_{i=1}^n\,\lb'_i\,p_{e'_i}$ 
with strictly positive coefficients $\lb_i,\lb'_i$, 
containing to index~$3$ the sublattice $\Lb_0$ generated 
by the $e_i$. Negating some $e_i$ if necessary, 
we may write $\Lb=\la\Lb_o,e\ra$ with $e=\frac{e_1+\dots+e_\ell}3$ 
for some $\ell\le n$. Then $\ell=n=6$, 
and the relation is proportional to Watson's relation 
$\sum_{i=1}^n\, p_{e_i}=\sum_{i=1}^n\,p_{e-e_i}$. 
\end{theorem}

\begin{proof} 
We have $\ell\ge 6$ by Lemma~\ref{lemwat} 
and $n=\ell$ by Corollary~\ref{corcyclic}. 

\vskip.1cm 

Next we bound from above the coefficients 
$A'_i=1-\sum_{j=1}^n\,(e'_i\cdot e_j^*)^2$ introduced 
in Lemma~\ref{lemAMB1}. 
We have $A'_i\le-1$ if $e'_i$ belongs to $\Lb_0$. Consider now a vector 
$e'_i\in\pm(e+\Lb_0)$, say, $e'_i=\frac{a_1 e_1+\dots+a_n e_n}3$ 
with $a_i\equiv 1\mod 3$. We have $A'_i=\frac{9-n}3$ if $e'_i=e$, 
$A'_i=\frac{6-n}9\le 0$ if $a_j=-2$ for exactly one index~$j$, 
i.e., if $e'_i=e-e_j$, and $A'_i<0$ otherwise. 
If none of the $e'_i$ is equal to $\pm e$, we have $A'_i\le 0$ 
for all~$i$, and even $A_i<0$ if $e'_i$ is not an $e-e_j$ or if $n>6$. 
Lemma~\ref{lemAMB1} then implies that $n=6$ and that all $e'_i$ 
are equal to some $e-e_j$. 
Our perfection relation is then of type Watson, as stated in the theorem. 

\vskip.1cm 

We must now consider the case when some $e'_i$ is equal to $\pm e$, 
say, $e'_1=e$, and prove that this is impossible. 
\vskip.1cm 

In the sequel, we normalize the coefficients $\lb_i,\lb'_i$ 
by the condition $\lb'_1=1$, and prove that the bound 
$\lb'_j A'_j\le-\frac 19$ holds for all $j\ge 2$. 
This will complete the proof of the theorem: indeed, we have
$A'_1=1-\frac n9$, 
hence $\sum_j\,\lb'_j A'_j\le(1-\frac n9)-\frac{n-1}9=\frac{10-2n}9<0$, 
which contradicts Lemma~\ref{lemAMB1}.

\vskip.1cm 

We first observe that for $j\ge 2$, $e'_j$ has at least two components 
different from $\pm\frac 13$. 
Otherwise, suppose that, say, $e'_2=e'_1+\lb e_1$. 
Then $e_1\cdot {e'_k}^*$ is zero for $k\ge 3$, 
hence $e'_k\in\Lb_0$ (because if $e'_k\in\pm e+\Lb_0$, 
then $e'_k\cdot e_1^*\equiv\pm\frac 13\mod\Z$). 
With the notation of Lemma~\ref{lemAMB4}, this implies $\nu_3=2$ 
and $\nu_1=n-2$, a contradiction. 

\vskip.1cm 

We now apply Lemma~\ref{lemAMB3},\,(4) with $k=1$: there exists $i_0$ 
such that $\lb'_j(e'_j\cdot e_{i_0}^*)^2\ge\frac 19$. Writing 
$A'_j=(1-\sum_{i\ne i_0}(e'_j\cdot e_i^*)^2)-(e'_j\cdot e_{i_0}^*)^2$, 
we see that the upper bound $\lb'_j A'_j\le-\frac 19$ 
holds whenever $1-\sum_{i\ne i_0}(e'_j\cdot e_i^*)^2\le 0$. 

This is clear if $e'_j$ belongs to $\Lb_0$ 
(there are at least two non-zero integral components). 

We now consider the case when $e'_j\in\pm e+\Lb_0$, 
say, $e'_j\in e+\Lb_0$. 
Its components are then $\frac 13,-\frac23,\frac 43,\dots$, 
and at least one for some $i\ne i_0$ differs from $\frac 13$, 
so that 
$1-\sum_{i\ne i_0}(e'_j\cdot e_i^*)^2\le 
1-\frac{4+(n-2)}9=\frac{7-n}9$, which suffices for $n\ne 6$. 

Let now $n=6$. The same conclusion clearly holds if at least 
two components with $i\ne i_0$ differ from $\frac 13$. 

Finally, in the remaining case, we use Lemma~\ref{lemAMB2},\,(4), 
taking $k=1$ and $i$ such that $e'_j\cdot e_i^*=\frac 13$. 
We obtain $\lb'_j\ge\frac 1{a_6}=\frac 13$. 
Since $A'_j\le 1-\frac{8+(n-2)}9=-\frac 13$, 
we again get the required upper bound $\lb'_j A'_j\le-\frac 19$. 
\end{proof}

\section{Lattices of index~4}\label{secindex4} 

We keep the notation of Section~\ref{seclatperfrel}, 
and now study perfection relations of the form 
$\sum_i\,\lb_i\,p_{e_i}=\sum_i\,\lb'_i\,p_{e'_i}$ 
with strictly positive coefficients $\lb_i,\lb'_i$ 
and system $\{e_i,e'_j\}$ of rank~$n$, 
when $[\Lb:\Lb_0]=4$ and $\Lb/\Lb_0$ is cyclic; 
see Section~\ref{seclat2elem} for non-cyclic $\Lb/\Lb_0$.
So we may write $n=m_1+m_2$ and $\Lb=\la\Lb_0,e\ra$ with 
{\small 
$$e=\frac{e_1+\dots+e_{m_1}+2(e_{m_1+1}+\dots+e_n)}4= 
\frac{f+e_{m_1+1}+\dots+e_n}2\,,$$ 
} 
\noi\hskip-.15cm 
where $f=\frac{e_1+\dots+e_{m_1}}2$. 
Note that the components of $e$ are non-zero by
Corollary~\ref{corcyclic} 
and recall that we have $n\ge 7$ and $m_2=1,2$ or $3$ if $n=7$, 
and of course $m_1\ge 4$. 

\smallskip 

We begin with three examples, which we shall prove 
to exhaust all possible perfection relations of the required form. 

\begin{example}\label{exn7d4a} 
{\small\rm 
Let $n=7$, $m_1=4$, and $m_2=3$. 
Set $e'_1=e$, $e'_2=e-e_6-e_7$, $e'_3=e-e_5-e_7$, $e'_4=e-e_5-e_6$, 
$e'_5=f-e_3-e_4$, $e'_6=f-e_2-e_4$, and $e'_7=f-e_2-e_3$. 
Then the vectors $e'_i$ are minimal 
and $\sum_i\,p_{e_i}=\sum_i\,p_{e'_i}$. 
} 
\end{example} 

\begin{proof} 
Making use of the $\D_4$-structures which exist in the spans 
of $e_1,e_2,e_3,e_4$ and of $f,e_5,e_6,e_7$, 
we obtain the two perfection relations 
$$\aligned
& p_{e_1}+p_{e_2}+p_{e_3}+p_{e_4}=p_f+p_{e'_5}+p_{e'_6}+p_{e'_7}\\ 
\text{ and }\\ 
& p_f+p_{e_5}+p_{e_6}+p_{e_7}=p_{e'_1}+p_{e'_2}+p_{e'_3}+p_{e'_4}\,. 
\endaligned$$ 
Eliminating $p_f$ proves what we want. 

\noi{\small 
[It is easily checked that the $14$ vectors $e_i,e'_i$ 
have perfection rank $r=13$. Note that the $15$ vectors 
$e_i,e'_i,f$ still have perfection rank $13$.] 
} 
\end{proof} 

\begin{example}\label{exn7d4b} 
{\small\rm 
Let $n=7$, $m_1=6$, and $m_2=1$. Watson's condition holds, 
and the corresponding perfection relation 
(cf. Proposition~\ref{propwat}) reads 
$$p_{e_1}+\dots+p_{e_6}+2 p_{e_7}= 
p_{e'_1}+\dots+p_{e'_6}+2 p_{e'_7}$$ 
with $e'_j=e-e_j$. 
} 
\end{example} 

\begin{example}\label{exn8d4} 
{\small\rm 
Let $n=8$, $m_1=8$, and $m_2=0$. Watson's condition holds, 
and the corresponding perfection relation 
(cf. Proposition~\ref{propwat}) reads 
$$p_{e_1}+\dots+p_{e_8}=p_{e'_1}+\dots+p_{e'_8}$$ 
with $e'_j=e-e_j$. 
} 
\end{example} 

\begin{theorem}\label{thindex4} 
Assume that $\Lb/\Lb_0$ is cyclic of order~$4$. 
Then $n=7$ 
\linebreak 
or $n=8$, 
and the perfection relation is one of the three relations 
\linebreak 
described 
in Examples \ref{exn7d4a}, \ref{exn7d4b} and \ref{exn8d4}. 
\end{theorem} 

For the proof, we shall have to consider several possibilities 
for the~$e'_j$. In particular, vectors of the set 
$$\cE=\big\{\pm\frac{e_1+\dots+e_{m_1} 
\pm 2 e_{m_1+1}\pm\dots\pm 2 e_n}4\big\}$$ 
will play an important r\^ole. For $x\in\Lb$, we denote 
by $\ord(x)$ the order ($1$, $2$ or $4$) of $x$ modulo $\Lb_0$. 
When $\ord(e'_j)\ne 1$, we assume that its component 
on $e_1$ is positive. 

\begin{lemma}\label{lemi4a} 
Let $i>m_1$. Then $e_i=\sum_{j=1}^n\,x_j e'_j$ satisfies 
Watson's equality $\sum_j\,\abs{x_j}=2$, with $x_j\ne 0$ 
\ifff $\ord(e'_j)=4$, and then $e'_j\cdot e_i^*=\pm\frac 12$. 
\end{lemma} 

\begin{proof} 
Recall (Lemma~\ref{lemAMB4}) that $N_4$ is the set of $e'_j$ 
with $\ord(r'_j)=4$. 
Set $S=\sum_j\,\abs{x_j}$ and $S_4=\sum_{e'_j\in N_4}\,\abs{x_j}$. 
By Lemma~\ref{lemAMB2},\,(2), we have 
$\sum_{j\mid e'_j\cdot e_i^*\ne 0}\,\abs{e'_j\cdot 
e_i^*}\,\abs{x_j}=1$, where $\abs{e'_j\cdot e_i^*}\ge\frac 12$ 
if $e'_j\in N_4$ and is $\ge 1$ otherwise. This implies 
$$\frac 12\,S_4+(S-S_4)\le 1\,,\eqno{(\ast)}$$ 
i.e. $S_4\ge 2(S-1)$, thus $S_4\ge 2$ since $S\ge 2$ 
by Watson's condition. Finally, equality holds in $(\ast)$, 
whence $\abs{e'_j\cdot e_i^*}=\frac 12$ if $e'_j\in N_4$, 
and $S=S_4$, whence $\abs{e'_j\cdot e_i^*}=0$ if $e'_j\notin N_4$, 
and eventually $S=2$. 
\end{proof} 


In order to use Lemma~\ref{lemAMB1}, we establish bounds 
for the $A'_j$ according to $\ord(e'_j)$: 

\smallskip 

\noi 
$A'_j\le -1$ if $\ord(e'_j)=1$, $A'_j\le 1-\frac{m_1}4\le 0$ 
if $\ord(e'_j)=2$, $A'_j=\frac{16-n-3 m_2}{16}$ if $e'_j\in\cE$, 
and $A'_j\le\frac{8-n-3 m_2}{16}$ if $\ord(e'_j)=4$ and $e'_j\notin\cE$, 
with equality \ifff $e'_j$ is of the form $x-e_i$ for some $i<m_1$ 
and $x\in\cE$. 
Lemma~\ref{lemAMB1} together with the data above shows 
that $m_2\le 3$ with equality only for $n=7$. 
We now study this case. 

\begin{lemma}\label{lemi4b} 
If $m_2\ge 3$, then $n=7$ and $m_2=3$, and the perfection relation 
is that of Example~\ref{exn7d4a}. 
\end{lemma} 

\begin{proof} 
Clearly, $m_2\ge 3$ implies that $n=7$ and $m_2=3$, 
and that $\max A'_j=0$, hence that all $A'_j$ are zero. 
This shows that either $e'_j\in\cE$, or $\ord(e'_j)=2$ 
and then $e'_j=\frac{e_1\pm e_2\pm e_3\pm e_4}2$. 
Since these last vectors have rank at most~$4$, 
there are at least three vectors $e'_j$ in $\cE$, 
and we may assume that $e'_1=e$. 
We now look at the other $e'_j$ lying in $\cE$. 

If, say, $e'_2=e'_1-e_7$, then $e_7=e'_1-e'_2$ has exactly 
two non-zero components on the basis $(e'_j)$. 
By Lemma~\ref{lemi4a}, $e'_1,e'_2$ are the only vectors 
of order~$4$, a contradiction. 

If, say, $e'_2=e'_1-e_5-e_6-e_7$, then a third vector $e'_3\in\cE$ 
would be of the form, say, $e'_3=e'_1-e_5-e_6=e'_2+e_7$, 
which contradicts the previous remark. 

Hence, two vectors $e'_j\in\cE$ must differ 
by exactly $2$ minus signs. Since no Watson relation 
as in Lemma~\ref{lemi4a} may involve exactly three vectors, 
there are exactly $4$~vectors $e'_j$ in $\cE$, say 
$e'_1=e$, $e'_2=e-e_5-e_6$, $e'_3=e-e_5-e_7$, and $e'_4=e-e_6-e_7$. 

The three vectors $e_5,e_6,e_7$ clearly belong to the span 
of $e'_1,e'_2,e'_3,e'_4$, and so does $f=2e'_1-e_5-e_6-e_7$ 
as well as the $f-e_i$ for $i\le 4$ (and $f-e_2-e_3-e_4=-(f-e_1)$). 

Hence the remaining three vectors $e'_j$ are those 
of Example~\ref{exn7d4a}. 
\end{proof}

\begin{lemma}\label{lemi4c} 
If no vector $e'_j$ belongs to $\cE$, then $n=8$ and the perfection 
relation is that of Example~\ref{exn8d4}. 
\end{lemma} 

\begin{proof} 
We know that if $e'_j$ has order $2$ (resp.~$1$) modulo $\Lb_0$, 
we have $A'_j\le 1-\frac{m_1}4\le 0$ (resp. $A'_j\le -1$); 
see Section~\ref{seclat2elem}. 

Let $e'_j$ be of order $4$ modulo $\Lb_0$, say, $e'_j\in e+\Lb_0$. 
since $e'_j$ does not belong to $\cE$, there exists $i\le m_1$ 
with at least one component different from $\frac 14$, 
thus equal to $-\frac 34,+\frac 54,\dots$, which implies 
$$A'_j\le 1-\frac{9+(m_1-1)+4 m_2}{16}=\frac{8-n-3 m_2}{16}\,.$$ 
At least one these $A'_j$ must be non-negative. 
We have $n\ge 7$ and $m_2\ge 1$ if $n=7$, 
which shows that we have indeed $n=8$ and $m_2=0$ and $A'_j=0$. 
Since $A'_j$ is then strictly negative if $e'_j$ 
has order $1$ or~$2$, all $e'_j$ must be of order~$4$ 
and all $A'_j$ must be zero. 
Hence all $e'_j$ are of the form $e-e_j$. 
\end{proof} 

\emph{From now on, we assume that $\cE$ contains some $e'_j$,
say $e'_1=e$, that $m_2\le 2$, that $\lb'_1=1$, 
and that $A'_1\ge 0$ (because $A'_1=\max A'_j$).} 

\begin{lemma}\label{lem1slash16} For every $e'_j\notin\cE$, we have $\lb'_j A'_j\le-\frac 1{16}$, 
with equality only if either $n=7$, $m_2=1$, $\nu_4=7$ 
and $e'_j$ is of the form $e'_1-e_i-e_7$ for some $i\le 6$, 
or perhaps $n=8$ and $m_2=1$. 
Moreover, if $n=8$ and $m_2=0$, we have the better inequality 
$\lb'_j A'_j\le-\frac 1{12}$. 
\end{lemma} 

\begin{proof} 
Let $e'_j=\sum_{i}\,x_i\,e_i$ not in $\cE$ (thus, $j\ge 2$). 
Suppose first that $(n,m_2)\ne(8,0)$. 
We use Lemma~\ref{lemAMB3} with $k=1$, 
setting $\s_i=\abs{e'_1\cdot e_i^*}\,\abs{e_i\cdot{e'_k}^*}$ 
and $X_i=\lb'_j\,\frac{x_i^2}{(e'_1\cdot e_i^*)^2}-1$. 
We have $\sum_i\,\s_i\,X_i=0$. The maximum of the $X_i$ 
is attained for some index $i_0\le m_1$ 
(if $i>m_1$, we have $x_i=0$ if $\ord(e'_j)\le 2$ 
and $x_i=\frac 12$ if $\ord(e'_j)=4$; see Lemma~\ref{lemi4a}). 
This maximum is non-negative, and even strictly positive 
since $e'_j\notin\cE$, hence $\lb'_j\,x_{i_0}^2>\frac 1{16}$. 

We now write 
$\lb'_j\,A'_j=\lb'_j(1-\sum_{i\ne i_0}\,x_i^2)-\lb'_j x_{i_0}^2$. 

If $\ord(e'_j)=1$, $1-\sum_{i\ne i_0}\,x_i^2$ is $\le 0$ 
since there are at least two non-zero, integral components. 
The same conclusion holds if $\ord(e'_j)=2$: we have 
$\sum_{i\ne i_0}\,x_i^2\ge\frac{m_1-1}4\ge 0$ 
(since $m_1=n-m_2\ge n-2\ge 5$). 
If $\ord(e'_j)=4$, and if at least two $x_i$, $i\le m_1$ 
are different from $\frac 14$, we have 
$1-\sum_{i\ne i_0}\,x_i^2\le\frac{9-n-3m_2}{16}\le 0$. 
If $x_{i_0}$ (say, $i_0=1$) is the unique component of $e'_j$ 
different from $\frac14$, say, 
$$e'_j=e'_1+(x_1-\frac 14)\,e_1-\sum_{i\in I}\,e_i\eqno(\dag)$$ 
with $I\subset\{m_1+1,\dots,n\}$. For all $i\ge 2$, 
$x_i=\pm e'_1\cdot e_i^*$, hence $X_i=\lb'_j-1$. 
Note that $I$ is not empty, and in particular that $m_2\ge 1$: 
otherwise, $e_1$ would have only 
two non-zero components and $e'_k$ would be in $\Lb_0$ 
for $k\ne 1,j$, whence $\nu_4=2$ and $\nu_1=n-2$, 
in contradiction with Lemma~\ref{lemAMB4}. 

The relation $\sum_i\,\s_i X_i=0$ (where $\sum_i\,\s_i=1$) 
now reads $\lb'_j=\dfrac 1{(16x_1^2-1)\s_1+1}$. 
We observe that for $i\in I$, $e_i$ does not belong to $\Lb'_0$: 
otherwise, by Lemma~\ref{lemi4a}, it would be of the form 
$e_i=\pm e'_k\pm e'_\ell$, where $e'_k$ and $e'_\ell$ are the only 
vectors of order~$4$ among the $e'_j$, which implies $e_i=e'_1-e'_j$, 
a contradiction. 
Hence the components of $e_i$ are bounded from above by $\frac 12$. 

From $(\dag)$, we obtain 
$(1-4x_1)\s_1=1-\sum_{i\in I}\,e_i\cdot{e'_1}^*$, 
where the right hand side is positive (because $m_2\le 2$) 
and indeed strictly positive (because $x_1\ne\frac14$) 
and at most~$1$ (because $e'_1\cdot e_i^*>0$). 
We have 
\linebreak 
$1-4x_i>0$, hence $1+4\a<0$, 
and the denominator of $\lb'_j$ is bounded from above 
by $-4x_1=4\abs{x_1}$, so that $\lb'_j\ge\frac 1{4\abs{x_1}}$. 
On the other hand, 
$A'_j=\frac{17-16x_1^2-n-3m_2}{16}$ is $\le 0$ 
(because $\abs{x_1}\ge\frac 34$), hence 
$$\lb'_j\,A'_j\le\frac{17-n-3m_2}{16\times 4\abs{x_1}}- 
\frac{4\abs{x_1}}{16}\,,$$ 
a decreasing function bounded from above by its value 
for $\abs{x_1}=\frac 34$, hence 
$\lb'_j\,A'_j\le\frac{8-n-3m_2}{48}\le-\frac 1{16}$, 
except if $n=7$ and $m_2=1$. 

In this case, we use a different argument. 
By \ref{lembases}, we have 
\linebreak 
$\lb'_j=\lb'_j\frac{\abs{e'_j\cdot e_7^*}}{\abs{e'_1\cdot e_7^*}}= 
\frac{\abs{e_7\cdot{e'_j}^*}}{\abs{e_7\cdot{e'_1}^*}}$, 
and by Lemma~\ref{lemi4a}, the right hand side is equal 
to $1$ if $e_7$ has at most $6$ non-zero components 
(in the subspace spanned by the corresponding $e'_j$, 
the index is at most~$3$), and to $1$ or $\frac 12$ 
if all components are non-zero, i.e. $\nu_4=7$. 
This implies $\lb'_j\ge\frac 12$ with equality only if $\nu_4=7$, 
and since $A'_j\le-\frac 2{16}$ with equality only 
for $\abs{x_1}=\frac34$, we obtain $\lb'_j\,A'_j\le-\frac 1{16}$, 
with equality as stated in the lemma. 

\smallskip 

We are left with the case when $(n,m_2)=(8,0)$. 
Then $e'_1=\frac{e_1+\dots+e_8}4$ is the only vector $e'_j$ 
belonging to~$\cE$. 
Let $j\ge 2$, and let $m$ be the minimum of the non-zero 
$\abs{e'_j\cdot e_i^*}$. By Lemma~\ref{lemAMB2}, 
we have $\lb'_j\ge\frac 1{24m}$. 
If $\ord{e'_j}=2\text{\,or\,}4$ (resp. $=1$), 
$e'_j$ has $8$ (resp. at least $2$) non-zero components, 
and thus $A'_j\le 1-8m^2$ (resp. $1-2m^2$), whence $A'_j\le 0$ 
since $e'_j\notin\cE$, and $\lb'j\,A'_j\le\frac{1-8m^2}{24m}$ 
(resp. $\frac{1-2m^2}{24m}$). 
The required inequality $\lb'_j\,A'_j\le-\frac 1{12}$ is proved when 
(1) $\ord(e'_j)=2$ because $m\ge\frac 12$; 
(2) $\ord(e'_j)=4$ and $m\ge\frac 34$; (3) $\ord(e'_j)=1$ and $m\ge 2$. 

If $\ord(e'_j)=4$ and $m=\frac 14$, then at least two components 
of $e'_j$ are distinct from $\frac 14$ 
(because $m_2=0$, see above, ``$I\ne\emptyset$''), 
which implies $A'_j\le-\frac 12$ and $\lb'_j\,A'_j\le-\frac 1{12}$. 

If $\ord(e'_j)=1$, $m=1$ and if $e'_j$ has at least three non-zero 
components, then $A'_j\le-2$ and the same conclusion again holds. 

There remains the case when, say, $j=2$ and $e'_2=e_1\pm e_2$, 
and thus $A'_2=-1$. We write 
$e_1=\frac{\sum_i\,a_k\,e'_k}d$ and $e_2=\frac{\sum_i\,a'_k\,e'_k}{d'}$ 
with coprime systems $\{d,a_k\}$ and $\{d',a'_k\}$. 
By Lemma~\ref{lemAMB2},\,(4), we have 
$$4\lb'_2=\frac{\abs{a_2}}{\abs{a_1}}=\frac{\abs{a'_2}}{\abs{a'_1}}\,.$$ 
We shall prove that $\frac{\abs{a_2}}{\abs{a_1}}\ge\frac 13$. 
This is true if $\abs{a_1}\le 3$, or $\abs{a_2}\ge 2$ 
(because we then have $\abs{a_1}\le\a_8=6$). 

We may now assume that $\abs{a_1}\ge 4$ $\abs{a_2}=1$, 
and note that the same hypotheses hold for $e'_2$. 
By Lemma~\ref{lemAMB2},\,(3) applied with $i=2$, 
we have $\frac 1d+\frac 1{d'}=1$, hence $d=d'=2$. 
By Lemma~\ref{lemAMB2},\,(2) applied with $j=1$, we have 
$\sum_{k,e'_k\cdot e_1^*\ne 0}\,\frac{\abs{a_k}}2\, 
\abs{e'_k\cdot e_1^*}=1$, and since $e'_k\cdot e_1^*\ge\frac 14$ 
and $\abs{a_1}\ge 4$, the previous formula gives 
$\frac\abs{a_1}4+\abs{a_2}+\sum_{k\ge 3}\,\frac{\abs{a_k}}4\le 2$. 
Since $d=2$ implies that $e_1$ has at least four non-zero components, 
this reads $\frac 52\le 2$. This contradicts the assumed values 
for $a_1$ and $a_2$, and completes the proof of the lemma. 
\end{proof} 

\smallskip 

\noi\emph{Proof of Theorem~\ref{thindex4}.} 

First consider the case when $e'_1$ is the only vector $e'_j$ in $\cE$. 
By Lemma~\ref{lem1slash16}, we have 
$\sum_j\,\lb'_j\,A'_j\le\frac{16-n-3m_2}{16}-\frac{n-1}{16}= 
\frac{17-2n-3m_2}{16}$. This is strictly negative, 
and thus contradicts Lemma~\ref{lemAMB1}, 
except if $(n,m_2)=(8,0)$ or $(7,1)$. 
In the first case, we may replace $\frac{n-1}{16}$ by $\frac{n-1}{12}$, 
which bounds the previous sum by $-\frac 1{12}<0$. 
If $(n,m_2)=(7,1)$, we have $\frac{17-2n-3m_2}{16}=0$, 
and by Lemma~\ref{lem1slash16}, this upper bound is strict 
unless all vectors $e'_j$, $j\ge 2$, have the form 
$e'_j=e'_1-e_i-e_7$, $i=1,\dots,6$. 
Replacing $e'_1$ by $e'_1-e_7$, we recover Example~\ref{exn7d4b}. 

\medskip 

From now on, we assume that $\cE$ contains at least 
two vectors $e'_j$, say, $e'_1=e$ and $e'_2$; 
this implies that $m_2\ge 1$, i.e. $m_2=1$ or~$2$. 
Recall that $\nu_d$ denotes the number of $e'_j$ 
of order $4$ modulo $\Lb_0$. 
We first show that there are actually exactly two (normalized) 
$e'_j$ in $\cE$, and that if $m_2=1$ (resp. $m_2=2$), 
then $e'_2=e'_1-e_n$ and $\nu_4=2$ 
(resp. $e'_2=e'_1-e_{n-1}-e_n$ and $\nu_4=4$). 

If $m_2=1$, then $e'_1-e_n$ is the unique possible choice 
for $e'_2$, and since then $e_n=e'_1-e'_2$, Lemma~\ref{lemi4a} 
shows $e'_1,e'_2$ are the only $e'_j$ of order~$4$. 

If $m_2=2$, suppose that (for instance) $e'_2=e'_1-e_n$. 
The argument above would imply that $e'_1,e'_2$ 
are the only $e'_j$ of order~$4$, and again by Lemma~\ref{lemi4a}, 
that $e_{n-1}$ would be of the form $e'_1\pm e'_2$, a contradiction. 
Hence $e'_2=e'_1-e_{n-1}-e_n$. By Lemma~\ref{lemAMB2},\,(1), 
the components of $e_{n-1}$ and $e_n$ on $e'_1$ (resp, $e'_2$) 
are positive (resp. negative), so that 
$\abs{e_{n-1}\cdot{e'_j}^*}+\abs{e_n\cdot{e'_j}^*}=1$ for $j=1,2$. 
Since Watson's equality is satisfied by $e_{n-1}$ and $e_n$, 
Lemma~\ref{lemwat} implies that $\abs{e_i\cdot{e'_j}^*}=\frac 12$ 
for $i=n-1,n$ and $j=1,2$, and that $e_{n-1}$ and $e_n$ 
are have denominator~$2$, hence have exactly $4$~non-zero components, 
i.e. $\nu_4=4$. 

When $m_2=2$, there are two vectors $e'_k$ of order $4$ outside $\cE$, 
which thus have at least one component distinct from $\frac 14$ 
on some $e_i$, $i\le m_1$, 
and then $A'_k\le 1-\frac{9+(m_1-1)+4m_2}{16}=\frac{2-n}{16}$. 
Lemma~\ref{lemAMB2},\,(4), applied with $h=1$ and $j=n$, 
reads $\lb'_k\,\frac{1/2}{1/2}=\frac{1/2}{1/2}$, i.e. $\lb'_k=1$. 
Hence, we have 
$\sum_j\,\lb'_j\,A'_j<2(\frac{10-n}{16})+2(\frac{2-n}{16})= 
\frac{12-2n}{16}<0$. This proves that $m_2=2$ is impossible. 

From now on, we assume that $m_2=1$. We have $e'_2=e'_1-e_n$ 
and $\abs{e_n\cdot{e'_1}^*}=1$. 
Fix $j\ge 3$ (thus, $\ord({e'_j}^*=1$ or~$2$. 
Lemma~\ref{lemAMB3} reads 
$$\sum\s_i\,X_i=0\ \text{ with }\ 
\s_i=\abs{e_i\cdot{e'_1}^*}\,\abs{e'_1\cdot e_i^*}\ \nd\ 
X_i=\lb'_j\,\frac{(e'_j\cdot e_i^*)^2}{(e'_1\cdot e_i^*)^2}-1\,.$$ 
Since $e_n\cdot{e'_j}^*=e'_j\cdot{e_n}^*=0$ (by Lemma~\ref{lemi4a}), 
we have $X_n=-1$. 
Since $\s_n=\frac12$ (and consequently $\sum_{i<n}\s_i=\frac 12$), 
we have $\frac 12=\sum_{i<n}\,\s_i\,X_i\le \frac 12\max_{i<n}\,X_i$. 
Hence there exists $i_0<n$ 
such that $\lb'_j(e'_j\cdot e_{i_0}^*)^2\ge\frac 18$. 
Since $1-\sum_{i\ne i_0}\,(e'_j\cdot e_{i_0}^*)^2\le 0$ 
(see the beginning of the proof of Lemma~\ref{lem1slash16}), 
we obtain $\lb'_j\,A'_j\le-\frac 18$, 
hence the upper bound 
$\sum_j\,\lb'_j\,A'_j\le 2\,\frac{13-n}{16}+(n-2)\,\frac{-1}8= 
\frac{15-2n}8$. This contradicts Lemma~\ref{lemAMB1}, 
except for $n=7$. 

There remains to deal with the case $n=7$ (and $m_2=1$). 

We again consider first a fixed $j\ge 3$, 
and set $\min_{i<n}\,\abs{e'_j\cdot e_i^*}$. 
By Lemma~\ref{lemAMB2},\,(4), we obtain 
$4x\lb'_j\ge\frac 1{\a_7}=\frac 14$, i.e. $\lb'_j\ge\frac 1{16x}$. 
If the absolute values of some component of $e'_j$ 
is $>x$, thus $\ge x+1$, then 
$A'_j\le 1-(5x^2+(x+1)^2)=-(6x^2+2x)$ if $\ord(e'_j)=2$, 
and $A'_j\le 1-(x^2+(x+1)^2)=-(2x^2+2x)$ if $\ord(e'_j)=1$, 
In both cases, we have $\lb'_j\,A'_j\le-\frac 14$. 
Otherwise, we have $A'_j=1-6x^2=-\frac 12$ if $\ord(e'_j)=2$ 
and $A'_j\le 1-2x^2=-1$ if $\ord(e'_j)=1$ 
(with equality only if $e'_j$ is of the form $e_i\pm e_k$). 
In the relation $\sum_{i\le 6}\,\s_i\,X_i=\s_7(=\frac 12)$, 
all $X_i$ with $i\le 6$ coincide with $16x^2\lb'_j-1$, 
whence $\lb'_j=\frac 1{8x^2}$ 
and $\lb'_j\,A'_j\le-\frac 14$ if $\ord(e'_j)=2$ 
and $\lb'_j\,A'_j\le-\frac 18$ if $\ord(e'_j)=1$. 
Since equality in Lemma~\ref{lemAMB4} does not hold 
(because $e'_1\cdot e_7^*=\frac 12\ne\pm\frac14$), 
we have $\nu_1<2\nu_4=4$, hence 
$\sum_j\,\lb'_j\,A'_j\le 2\,\frac 6{16}-3\,\frac 18-2\,\frac 14= 
-\frac 18<0$, which again contradicts Lemma~\ref{lemAMB1}. 
This completes the proof of Theorem~\ref{thindex4}. 
\qed 

\begin{corol}\label{cor7} 
Let $\Lb$ be a lattice of dimension $n\le 7$ endowed with a perfection 
relation $\sum_{i=1}^n\,\lb_i\,p_{e_i}=\sum_{i=1}^n\,\lb'_i\,p_{e'_i}$ 
with strictly positive coefficients $\lb_i,\lb'_i$, 
where $\{e_i,e'_i\}$ is a rank~$n$ set of minimal vectors of~$\Lb$. 
Let $\Lb_0$ and $\Lb'_0$ be the lattices generated by 
the $e_i$ and the $e'_i$, and assume that $\Lb_0+\Lb'_0=\Lb$. 
Then $\Lb$ is similar to one of the lattices $\D_4$, $\D_6$ or $\E_7$, 
endowed with a perfection relation as in Theorem~\ref{th2elem}, 
or $n=6$ (resp. $n=7$), and $\Lb$ is as described 
in Theorem~\ref{thindex3} (resp. Theorem~\ref{thindex4}). \qed 
\end{corol}

\section{Complements}\label{seccompl} 

\subsection{Cells}\label{subseccell} 

The space of positive definite quadratic forms 
has a natural structure of an (infinite) cell complex. 
In the dictionary lattices --- quadratic forms, the set of cells 
which are equivalent to a given one corresponds to a \emph{minimal class}; 
see \cite{M}, Section~1.7 and Chapter~9. This dictionary relies 
on the choice of a basis $\cB$ for $E$, with which we attach to any 
$x\in E$ the column $X$ of its components in~$\cB$. 
For a unitary vector~$x$, we then have $\Mat(p_x,\cB^*,\cB)=X\,\ta X$. 
Perfection relations then read $\sum_x\,\lb_x\,X\,\ta X=0$. 
This shows that the space of perfection relations on a lattice $\Lb$ 
is an invariant of the class $\cC$ of~$\Lb$. 

\smallskip 

The dimension of a cell 
is its \emph{perfection co-rank}, equal to $\frac{n(n+1)}2-r$, 
where $r$ is 
the perfection rank of any lattice in the class. 
An inclusion $\cC\supset\cC'$ between cells is equivalent to the opposite 
inclusion $S\subset S'$ on sets of minimal vectors, and induces 
on the set of (minimal) classes an ordering relation denoted 
by ``$\prec$'', for which the maximal classes are the perfect ones, 
corresponding to cells of dimension~$0$. 
The space of perfection relations of a class $\cC$ embeds canonically 
in the corresponding set of any class $\cC'\succ\cC$. 
In particular, all perfection relations in a given dimension 
can be constructed using only perfect lattices having this dimension 
(and of course $s>r=\frac{n(n+1)}2$). 
In dimensions $n\le 6$, these lattices are $\D_4$ and $\D_5$ 
($\imath=2$), $\E_6$, $\E_6^*$ and the perfect, non-extreme 
$6$-dimensional lattice $P_6^4$ ($\imath=3$), 
and $\D_6$ ($\imath=2^2$).

\smallskip 

Given an integer $d\ge 1$ and a code $C$ of length $n$ over $\Z/d\Z$, 
we say that a pair $(\Lb,\Lb')$ of lattices is \emph{admissible for $C$} 
if $\Lb$ is well rounded of dimension~$n$, $\Lb'$ is a sublattice 
of $\Lb$ generated by minimal vectors of $\Lb$, 
and there exists a basis for $\Lb$ inducing an isomorphism 
$\Lb/\Lb'\simeq C$; 
the list of $\Z/d\Z$-codes possessing admissible pairs $(\Lb,\Lb')$ 
is displayed in Table~11.1 of \cite{M1} up to length~$8$. 

This notion of admissible pairs is again a class invariant, 
in the sense that if $(\Lb,\Lb')$ is admissible for $C$, 
then every lattice $L$ in the class $\cC$ of $\Lb$ 
contains a sublattice~$L'$ such that $L/L'$ defines the same code 
(and this result even holds for every class $\cC'\succ\cC$). 
Moreover, the averaging argument developed in Section~8 of~\cite{M1} 
shows that the set of cells admissible for $C$, if non-empty, 
is the set of cells $\cC\succ\cC_0$ for a uniquely defined 
cell~$\cC_0$. 

\smallskip 

\noi{\bf Example.}{\small 
When Watson's equality holds, Lemma~\ref{lemwat} 
shows that $\Lb$ must contain at least $2\ell$~vectors 
and even $3\ell$ if $a_\ell=\pm\frac d2$. 
When $d\le 5$, $\cC_0$ exists and its minimal vectors are exactly 
those listed Lemma~\ref{lemwat}. However, for $d=6$ and $n=8$ 
(the smallest dimension for which index~$6$ may occur), 
there is a unique class satisfying Watson's equality, 
and this is that of $\E_8$, with $s=120$, much larger that $3n=24$; 
see Subsection~\ref{subsecdim8} below.} 

\smallskip 

\noi{\bf Remark.}{\small\emph{ 
All the perfection relations that we have classified 
in the previous sections involve only vectors of $\cC_0$}. 
We do not know whether this is general.} 

\subsection{A glance at dimension 8}\label{subsecdim8} 

We now consider perfection relations 
$\sum_{i=1}^n\,\lb_i\,p_{e_i}= \sum_{i=1}^n\,\lb'_i\,p_{e'_i}$ 
with strictly positive coefficients $\lb_i,\lb'_i$ 
and independent unitary vectors $e_1,\dots,e_n$ 
(and thus also $e'_1,\dots,e'_n$) 
in the case when~$n=8$; the notation 
$\Lb_0$, $\Lb'_0$, $\Lb$ is as usual. 

It results from the previous sections that $\imath=[\Lb:\Lb_0]\ge 4$, 
and that if $\Lb/\Lb_0$ is cyclic of order~$4$, 
then $\Lb$ is as in Example~\ref{exn8d4}. 
We also know that if $\Lb/\Lb_0$ is $2$-elementary, 
then $\Lb$ is similar to a direct sum $\D_4\bgop\D_4$ 
(the perf-reducible case) if $\imath=4$, to $\D_8$ if $\imath=8$, 
and to $\E_8$ if $\imath=16$ (one regular and one irregular relation). 

For the other $8$-dimensional perfection relations, $\Lb/\Lb_0$ 
must be of type $(5)$, $(6)$, $(4,2)$ or $(3,3)$; 
we give some examples below. 

\smallskip 

$\bullet$ $[\Lb:\Lb_0]=5$. Two relations are mentioned 
in Section~\ref{seclatperfrel}, 
with $(m_1,m_2)=(6,2)$ (a relation of type Watson) 
and $(4,4)$ (related to an identity of Zahareva). 
We conjecture that these are the only examples, 
and even more generally, that in larger dimensions, 
there are only two other cases, namely the relations 
of type Watson for $(m_1,m_2)=(8,1)$ and $(10,0)$. 

\smallskip 

$\bullet$ $[\Lb:\Lb_0]=6$. There is one relation of Watson type 
corresponding to $(m_1,m_2,m_3)=(5,2,1)$. 
In this case (as for the quotients of type $(3,3)$), 
the corresponding lattices are similar to $\E_8$, 
and thus a lot of relations could occur (they span a space 
of dimension $120-36=84$). 

Another regular relation exists in case $(m_1,m_2,m_3)=(2,4,2)$. 
The lattice $\Lb$ is generated over $\Lb_0$ 
by $e=\frac{e_1+e_2+2(e_3+e_4+e_5+e_6)+3(e_7+e_8)}6$. 
Besides the $8$ vectors $e_i$, the Watson equalities produce 
the $8$~vectors $\frac{e_1\pm e_2\pm e_7\pm e_8}2$ 
and the $6$~vectors $e'-e_1$, $e'-e_2$ and $e'_i=e'+e_i$, 
$i=3,4,5,6$, where $e'=\frac{e_1+e_2-e_3-e_4-e_5-e_6}3$, 
and besides these $8+8+6=22$ vectors, $\Lb$ contains 
$6$~extra minimal vectors, 
namely $e'$, $e'-e_1-e_2$, $e'_7=e$, $e'_8=e-e_7-e_8$, 
$e-e_7$ and $e-e_8$. Setting moreover 
$e'_1=\frac{(e_1-e_2)+(e_7-e_8)}2$ and 
$e'_2=\frac{(e_1-e_2)-(e_7-e_8)}2$, we have the relation 
$\sum_{i=1}^8\,p_{e_i}=\sum_{i=1}^8\,p_{e'_i}$, 
which again exists in~$\cC_0$. 

\smallskip 

$\bullet$ $[\Lb:\Lb_0]=8$ (quotient of type $(4,2)$). 
We do not know any example. There are three codes over $\Z/4\Z$, 
and we can only prove that the first which appears in \cite{M1}, 
Table~11.1 is impossible. 
Indeed, we then have $\Lb=\la\Lb_0,e,f\ra$ with 
$$e=\frac{e_1+e_2+e_3+e_4+2(e_5+e_6+e_7)}4\ \nd\ 
f=\frac{e_1+e_2+e_5+e_8}2\,,$$ 
and we observe that for any possible $e'_j$, we have $A'_j\le 0$, 
with equality only if $e'_j$ is one of the $48$ pairs of vectors 
of $\cC_0$ listed in \cite{M1}, Section~10. 
We can then conclude using the \emph{PARI} package. 

\smallskip 

$\bullet$ $[\Lb:\Lb_0]=9$ (quotient of type $(3,3)$). 
The irregular relation on $\E_8$ produces an example. 
We do not know whether other relations exist. 

\subsection{Group actions}\label{subsecgroups} 

Given a lattice $\Lb$, any subgroup $G$ of $\Aut(\Lb)$ 
acts of the real vector spaces $\cP$, the span in $\End^s(E)$ 
of the projections $p_x$, $x\in S(\Lb)$ and $\cRR$ of perfection 
relations of $\Lb$. When this last action is irreducible, 
any non-trivial perfection relation together with its conjugates 
under $G$ spans $\cRR$. 
Note that the symmetric square $\Sym^2(E)$ of the representation 
afforded by the action of $G$ on~$E$ is generated by the set 
of projections to the lines of~$E$ and thus coincides with $\cP$ 
when $\Lb$ is perfect. 

\smallskip 

These remarks apply to the perfect lattices 
$\E_6$, $\E_6^*$, $\E_7$, $\E_8$, for which $\cRR$ has dimensions 
$15$, $9$, $35$ and $84$ respectively, taking for $G$ 
the corresponding Weyl groups. 
Reduction modulo~$2$ gives isomorphisms $W(\E_6)\simeq\Oo_6^-(2)$, 
$W(\E_7)\simeq\Oo_7(2)$ and $W(\E_8)\simeq\Oo_8^+(2)$. 
An inspection of the character tables displayed in \cite{ATLAS} 
(in the unitary notation $U_4(2)$ for $W(\E_6)$) 
proves that $\cRR$ is irreducible in all cases. 

\smallskip 

For $\Lb=\D_n\subset\Z^n$ ($n\ge 4$), 
we take for $G$ the symmetric group $S_n$. 
Let $p_n$ be the permutation representation afforded by the canonical 
basis for $\Z^n$. Write $p_n=1+u_n$. Using \cite{F-H}, exercise~4.19, 
we check that $\Sym^2(u_n)$ is of the form $p_n+q_n$, 
where $q_n$, of dimension $\frac{n(n-1)}2-n=\frac{n(n-3)}2$, 
is the irreducible representation attached to the partition 
$[n-2,2]$ of~$n$. 

\smallskip 

As a consequence, we see that for the lattices considered above 
except $\E_6^*$, the space of perfection relations is generated 
by the conjugates of one relation coming from one cross-section 
$\D_4$. For $\E_6$, $\E_6^*$, $\E_7$, $\E_8$, we may take 
a Watson relation coming from a $6$-dimensional section. 

\subsection{Some other relations from dimension 8}\label{subsecotherrel} 

Using {\em PARI}, we have made an exploration of the $8$-dimensional 
perfect lattices with $s=37$ (thus, $s=r+1$). 
There are $2033$ such lattices, for which there exists a unique 
(up to proportionality) perfection relation, that we write in the form 
$\sum_{i=1}^{n1}\,\lb_i\,p_{e_i}=\sum_{j=1}^{n2}\,\lb'_j\,p_{e'_j}$ 
with strictly positive $\lb_i,\lb'_j$. For three of them, 
we have $\{n_1,n_2\}=\{14,10\}$. For the remaining $2030$ relations, 
we have $n_2=n_1$, with the following numbers of occurrences: 
$n_1=6:1450$; $n_1=8:404$; $n_1=9:56$; $n_1=10:87$; 
$n_1=11:11$; $n_1=12:10$; $n_1=13:12$. 
Relations with $n_1=6$ come from a Watson's relation with index~$3$. 
Those with $n_1=8$ share out among the following three types: 

\vskip.1cm 

{\small 
\ctl 
{index~$4$\,: $(m_1,m_2)=(8,0)$ ($64$);} 

\vskip.1cm 

\ctl{ 
index~$5$\,: $(m_1,m_2)=(4,4)$ ($338$); $(m_1,m_2)=(6,2)$ ($2$).}
\noi[That no other type with $n_1\le 8$ may occur 
can be proved {\sl a priori\/}, using the equality $s=r+1$.]} 

\smallskip 

Among the relations with $n_1=9$, 2 arise from $6$-dimensional 
lattices, all with $s=18$ (thus, the relation exists 
once more in $\cC_0$). This class has index~$2$, and we indeed have 
$[\Lb:\Lb_0]=[\Lb:\Lb'_0]=2$. Such a relation is easily constructed 
in both the perfect lattices $\E_6$ and $\E_6^*$, which needs 
the consideration of index~$3$, \emph{but not in~$\D_6$}. 
(For the lattices $\D_n$, $n_1$ and $n_2$ must be even.) 
Thus the relation cannot be constructed within $\cC_0$ 
using relations described in the previous sections. 

\smallskip 

Among the relations with $n_1=10$, $13$ out of $87$ come from a 
\linebreak 
$7$-dimensional lattice. 
In all other cases, dimension~$8$ is needed.

\vskip.6cm 

\vbox{{\small 
\noi 
Anne-Marie {\sc Berg\'e} {\&} Jacques {\sc Martinet} , 
A2X, Inst. Math., Universit$\acute{\rm e}$ Bordeaux 1, 
351, cours de la Lib$\acute{\rm e}$ration 
33405 Talence cedex, France. 

\noi 
{\it E-mails:\/} 
{\tiny\textsf 
berge@math.u-bordeaux.fr\,; martinet@math.u-bordeaux.fr\,.}  
}} 
\vfil\eject 

\end{document}